\documentclass[11pt]{amsart}
\usepackage[a4paper,margin=25mm]{geometry}
\usepackage{amsmath,amssymb,amsthm,booktabs,mathtools}
\usepackage[colorlinks=true,linkcolor=blue,citecolor=blue,urlcolor=blue]{hyperref}
\allowdisplaybreaks[2]
\numberwithin{equation}{section}
\newtheorem{theorem}{Theorem}[section]
\newtheorem{lemma}[theorem]{Lemma}
\newtheorem{proposition}[theorem]{Proposition}
\newtheorem{remark}[theorem]{Remark}
\newcommand{\F}{\mathbb F}
\newcommand{\Sp}{\operatorname{Sp}}
\newcommand{\Sw}{\operatorname{Sw}}
\newcommand{\Tr}{\operatorname{Tr}}
\title[Self-adjoint skew quartics]{Self-adjoint skew quartics, symplectic monodromy, and zeta functions}
\hypersetup{pdftitle={Self-adjoint skew quartics, symplectic monodromy, and zeta functions},
  pdfauthor={Tomoya Hayata and Takahiro Tsushima}}
\author{Tomoya Hayata}
\address[Tomoya Hayata]{Keio University School of Medicine,
4-1-1 Hiyoshi, Kohoku-ku, Yokohama, 223-8521, Japan}
\address[Tomoya Hayata]{RIKEN Center for Interdisciplinary Theoretical
and Mathematical Sciences (iTHEMS), RIKEN, Wako 351-0198, Japan}
\email{hayata@keio.jp}

\author{Takahiro Tsushima}
\address[Takahiro Tsushima]{Keio University School of Medicine,
4-1-1 Hiyoshi, Kohoku-ku, Yokohama, 223-8521, Japan}
\email{tsushima@keio.jp}
\makeatletter
\@ifundefined{subjclassname@2020}{%
  \@namedef{subjclassname@2020}{2020 Mathematics Subject Classification}}{}
\makeatother
\subjclass[2020]{Primary 11T06; Secondary 16S36, 11G20, 20C33}
\keywords{Finite fields, skew polynomials, irreducibility, local Fourier transform,
Deligne--Lusztig characters, zeta functions}
\date{}
\begin{document}
\begin{abstract}
For an odd prime power $P$, let $A_H(P)$ be the number of
$a\in\F_{P^H}$ for which the self-adjoint skew quartic
$t^4+t^3+at^2+t+1$ is irreducible in
$\F_{P^H}[t;\sigma]$, with $\sigma(c)=c^P$.
We determine $A_H(3)$ for every $H\ge1$ and prove that its
minimal constant-coefficient recurrence has order $24$. Local
ramification gives the degrees of the Steinberg and regular
cuspidal Frobenius polynomials. Duality and finitely many exact
counts determine both polynomials.

For $P=p^f$ with $p\nmid f$, we construct an explicit
genus-$P(P+1)/2$ flag curve whose zeta numerator is the Steinberg
Frobenius polynomial. At $P=3$ this identifies the degree-$12$
factor computed above with the zeta numerator of a genus-$6$
curve. For general $P$ the flag curve and the regular cuspidal
factors give a rational expression for the counting series. We
evaluate the linear coefficients of the two factors by an
elliptic-curve point count and express their quadratic
coefficients through a single finite-field character sum.
\end{abstract}
\maketitle

\section{Introduction}\label{sec:intro}
Let $P$ be an odd prime power and set $Q=P^H$. In the skew polynomial
ring $\F_Q[t;\sigma]$, where $tc=c^Pt$, consider the one-parameter
self-adjoint family
\[
 f_a(t)=t^4+t^3+at^2+t+1,\qquad a\in\F_Q.
\]
The Laurent polynomial $f_at^{-2}$ is fixed by the anti-involution
that fixes coefficients and sends $t$ to $t^{-1}$. We write $A_H(P)$
for the number of parameters for which $f_a$ is irreducible. Our aim
is to determine these counts over extensions of the parameter field
and to identify the geometric and representation-theoretic structures
responsible for them.

The starting point is that the self-adjoint coefficient pattern gives
the roots of the associated linearized polynomial a nondegenerate
alternating pairing. Irreducibility in the Ore ring is then equivalent
to irreducibility of Frobenius on this four-dimensional symplectic
root space. Thus the counting problem selects particular Frobenius
classes in a finite symplectic group. Expanding their indicator in
irreducible characters turns the count over each extension field
into a sum of cohomological traces.

The first result is a complete evaluation at $P=3$.
Theorem~\ref{thm:main} gives $A_H(3)$ for all $H$ in terms of
explicit polynomials $N(u)$ and $D(u)$ of degrees $12$ and $10$.
The former is the Steinberg factor and the latter is the regular
cuspidal factor. We compute the ramification of the root local
system at infinity using the Weil local systems of Katz--Tiep
\cite{KT} and local Fourier transformation
\cite{Laumon,AS}. The resulting conductor degrees and Frobenius
duality give a finite system of equations for the coefficients.
The counts through $H=11$ determine both factors, and the count
at $H=12$ checks the reconstruction. Proposition~\ref{prop:no-cancel-three}
then proves that the minimal recurrence has order $24$.

The second result gives a geometric interpretation of the
Steinberg factor for $P=p^f$ with $p\nmid f$. Theorem~\ref{thm:flagcurve}
constructs a smooth projective curve $\mathcal B_P$ of genus
$P(P+1)/2$ from flags consisting of a root line in a Lagrangian
root plane. Its zeta numerator is the Steinberg Frobenius
polynomial $S_P(u)$. In particular, $S_3(u)=N(u)$ identifies
the numerator factor already computed in Part I with the
zeta numerator of $\mathcal B_3$. The construction of the flag
curve uses the same monodromy and ramification calculation as
Part I.

With $m=(P^2+1)/2$ and $\kappa=(P^2-1)/8$, the general counting
formula is
\[
 m A_H(P)=\kappa\#\mathcal B_P(\F_{P^H})
                 +\sum_{\lambda}\lambda^H,
\]
where $\lambda$ runs through the inverse roots of the regular
cuspidal factor $C_P(u)$, with multiplicity. The degrees of
$S_P$ and $C_P$ are $P^2+P$ and $\kappa(P^2+P-2)$,
respectively. Part II gives explicit formulas for the first two
coefficients of each factor, including in characteristics $3$
and $5$, and relates the full counting series to a quotient of
covering-curve zeta functions.

Ore's theory \cite{Ore} and the reduced-norm methods of
Caruso--Le Borgne \cite{CL} relate skew-polynomial factorization
to semilinear algebra. The self-adjoint trinomials
$x^{2n}+ax^n+1$ studied by Hayata and Tsushima
\cite{HTtrinomials} give rise to quadratic reduced norms.
For the present quartic, the symplectic character expansion
\cite{DL} instead isolates a Steinberg term and regular cuspidal
terms. The appendix records the alternating product of flag-cover
zeta functions associated with the Steinberg representation in
higher rank.

\part{Exact enumeration over extensions of \texorpdfstring{$\F_3$}{F3}}
\section{Companion modules and the formula over
  \texorpdfstring{$\F_3$}{F3}}\label{sec:setup}
Let $P$ be an odd prime power, put $Q=P^H$, and let
$R=\F_Q[t;\sigma]$, where $\sigma(c)=c^P$ and
\[
 tc=c^Pt\qquad(c\in\F_Q).
\]
We write $A_H(P)$ for the number of parameters $a\in\F_Q$ for which
\[
 f_a(t)=t^4+t^3+at^2+t+1
\]
is an atom of $R$. We first work with an arbitrary odd $P$, since the
character identity underlying the enumeration is uniform in $P$;
we then specialize to $P=3$.
The units of $R$ are exactly the nonzero constants. A nonconstant
polynomial is an \emph{atom} if it cannot be written as a product of
two nonunits. This is irreducibility in the Ore ring, not in the
commutative polynomial ring $\F_Q[t]$.

\subsection*{Companion modules and the root space}
The quotient $R/Rf_a$, with $Rf_a$ the left ideal generated by $f_a$,
is a left $R$-module of dimension four over $\F_Q$.
Left multiplication by $t$ is $\sigma$-semilinear; in the basis
$1,t,t^2,t^3$ its companion matrix is
\[
 C_a=\begin{pmatrix}
 0&0&0&-1\\ 1&0&0&-1\\ 0&1&0&-a\\ 0&0&1&-1
 \end{pmatrix}.
\]
Right factors of $f_a$ correspond to submodules of this cyclic left
module. Consequently, $f_a$ is an atom if and only if $R/Rf_a$ is
simple.

Consider also the linearized polynomial and its root space
\[
 L_a(X)=X^{P^4}+X^{P^3}+aX^{P^2}+X^P+X,
 \qquad V_a=\ker(L_a)\subset\overline{\F}_P.
\]
Because the coefficient of $X$ is nonzero, $L_a$ is separable,
and $V_a$ is a four-dimensional $\F_P$-vector space.
The linearization map $\sum_i c_i t^i\mapsto\sum_i c_i X^{P^i}$
takes skew multiplication to composition. A monic proper right
factor $g$ of $f_a$ has nonzero constant coefficient, so its
linearization has a nonzero proper Frobenius-stable root space in
$V_a$. Conversely, let $W'\subset V_a$ be a nonzero proper
$F_Q$-stable $\F_P$-subspace. Its monic root polynomial
$\prod_{w\in W'}(X-w)$ is linearized. Indeed, starting with $X$,
adjoining a basis vector $h$ to a subspace with root polynomial
$D_{W'}(X)$ replaces that polynomial by
$D_{W'}(X)^P-D_{W'}(h)^{P-1}D_{W'}(X)$.
Frobenius permutes its roots, so its coefficients lie in $\F_Q$.
Let $g$ be the associated monic skew polynomial. Right division
$f_a=qg+r$ gives $\deg r<\deg g$. The linearization of $r$
vanishes on $W'$, which has $P^{\deg g}$ elements. Its degree
if nonzero is at most $P^{\deg r}$, so $r=0$. Thus $f_a$ is
an atom exactly when $F_Q$ has no nonzero proper invariant
subspace in $V_a$. For a linear operator this is equivalent to
having irreducible characteristic polynomial, by viewing the
space as an $\F_P[Z]$-module with $Z$ acting by $F_Q$. Therefore
\begin{equation}\label{eq:atom-frob-criterion}
 f_a\text{ is an atom}
 \quad\Longleftrightarrow\quad
 \det(Z-F_Q\mid V_a)\text{ is irreducible over }\F_P.
\end{equation}
Here and below, ``irreducible'' for the Frobenius action means
irreducibility of its four-dimensional characteristic polynomial.

For clarity, put
$\mathbf v(x)=(x,x^P,x^{P^2},x^{P^3})^{\mathsf t}$.
The equation $L_a(x)=0$ gives
\[
 \mathbf v(x^P)=C_a^{\mathsf t}\mathbf v(x).
\]
Successive application of Frobenius therefore yields
\[
 \mathbf v(x^Q)=M_{a,H}^{\mathsf t}\mathbf v(x),\qquad
 M_{a,H}=C_aC_{a^P}\cdots C_{a^{P^{H-1}}}.
\]
If $x_1,\ldots,x_4$ is an $\F_P$-basis of $V_a$, the matrix with
columns $\mathbf v(x_i)$ is a nonsingular Moore matrix. The preceding
identity shows that the matrix of $F_Q$ on $V_a$ is similar over
$\overline{\F}_P$ to $M_{a,H}^{\mathsf t}$.
In particular, the two operators have the same characteristic
polynomial in $\F_P[Z]$. This gives the concrete companion-matrix test for
\eqref{eq:atom-frob-criterion}.

\subsection*{The symplectic structure}
A direct computation gives
\[
 C_a^{\mathsf t}JC_a=J,\qquad
 J=\begin{pmatrix}
 0&0&1&-1\\0&0&0&1\\-1&0&0&0\\1&-1&0&0
 \end{pmatrix},\qquad \det J=1.
\]
This is a common alternating form for the companion matrices.
For the root-space realization, it is convenient to use the
following (equivalent) symplectic structure:
\begin{align*}
 B(x,y)={}&x^{P^2}y-xy^{P^2}
       +x^{P^2}y^P-x^Py^{P^2}
       +x^{P^3}y^P-x^Py^{P^3}.
\end{align*}
Its defining identity is
\begin{equation}\label{eq:B-artin-schreier}
 B(x,y)^P-B(x,y)
 =L_a(x)y^{P^2}-x^{P^2}L_a(y).
\end{equation}
Thus $B(x,y)\in\F_P$ for $x,y\in V_a$.
In Moore coordinates it is represented by the alternating matrix
\[
 K=\begin{pmatrix}
 0&0&-1&0\\0&0&-1&-1\\1&1&0&0\\0&1&0&0
 \end{pmatrix},\qquad \det K=1.
\]
Since the Moore matrix of any $\F_P$-basis of $V_a$ is invertible,
$B$ is nondegenerate on $V_a$. Furthermore,
$B(x^Q,y^Q)=B(x,y)^Q=B(x,y)$ for $x,y\in V_a$.
Hence the arithmetic Frobenius acts through
$\Sp(V_a,B)\simeq\Sp_4(\F_P)$.

As $a$ varies, the root spaces form a rank-four \'etale
$\F_P$-local system $\mathcal V$ on $\mathbb A^1_{\F_P}$,
and $B$ defines a perfect alternating pairing on it.
Fix a geometric generic point $\bar\eta$ and a symplectic basis
of $\mathcal V_{\bar\eta}$. We obtain a monodromy representation
\[
 \rho:\pi_1(\mathbb A^1_{\F_P},\bar\eta)
 \longrightarrow G:=\Sp_4(\F_P).
\]
Changing the basis conjugates $\rho$, so its Frobenius conjugacy
classes and all the character values below are well defined.
For $g\in G$ define $1_{\mathrm{irr}}(g)=1$ when its natural
four-dimensional representation is irreducible, and $0$ otherwise.
By \eqref{eq:atom-frob-criterion},
\begin{equation}\label{eq:A-as-frob-sum}
 A_H(P)=\sum_{a\in\F_{P^H}}
 1_{\mathrm{irr}}\bigl(\rho(\operatorname{Frob}_a)\bigr),
\end{equation}
with the arithmetic Frobenius convention.
Since $g$ is irreducible if and only if $-g$ is irreducible,
$1_{\mathrm{irr}}$ descends to
$S=G/\{\pm I\}=\operatorname{PSp}_4(\F_P)$.

\subsection{The irreducibility indicator and a Coxeter-torus identity}
Let $T\subset G$ be a Coxeter torus. Concretely, on the natural
four-dimensional module it is realized by the norm-one subgroup
\[
 T\simeq\{z\in\F_{P^4}^{\times}:z^{P^2+1}=1\},
 \qquad |T|=P^2+1.
\]
The Frobenius orbit of a noncentral element of $T$ has length four:
indeed, if $z\in T$ and $z^{P^2}=z^{-1}$, then an orbit of length
at most two forces $z=\pm1$. Thus every $t\in T\setminus\{\pm I\}$
is regular semisimple and acts irreducibly on the natural module.
Conversely, any irreducible $g\in G$ is regular semisimple and belongs
to a conjugate of $T$: its eigenvalues are
$z,z^P,z^{P^2},z^{P^3}$, and symplectic self-duality forces
$z^{P^2}=z^{-1}$. To see the conjugacy assertion, identify the
natural module with $\F_{P^4}$ so that $g$ acts by multiplication
by $z$. The adjoint on $\F_P[g]=\F_{P^4}$ is
$\bar x=x^{P^2}$ because it sends $z$ to $z^{-1}$.
By nondegeneracy of the trace pairing, every compatible
alternating form is $\Tr_{\F_{P^4}/\F_P}(\lambda x\bar y)$ with
$\bar\lambda=-\lambda$ and $\lambda\ne0$. The ratio of any two
such $\lambda$ belongs to $\F_{P^2}^{\times}$, and multiplication
by an element of $\F_{P^4}^{\times}$ of that norm transports one
form to the other. This gives a symplectic conjugacy of the
norm-one tori.

Put
\[
 T_0=T/\{\pm I\}\subset S,\qquad
 m=|T_0|=\frac{P^2+1}{2},\qquad
 \kappa=\frac{m-1}{4}=\frac{P^2-1}{8}.
\]
As $P$ is odd,
$\gcd(m,P^2-1)=1$. The rational Weyl group of $T$ is the finite quotient
$N_G(T)/T$, where $N_G(T)$ is the normalizer of $T$ in $G$.
For this Coxeter torus it is cyclic of order four; its action
on $T$ and on the character group $\widehat T$ is induced by
$P$-power Frobenius.
Every nontrivial character of $T_0$ therefore has an orbit of
length four. Let $\mathcal O$ be the set of these orbits;
then $|\mathcal O|=\kappa$.
For each orbit choose a representative $\theta$ and inflate it to
$T$, so that $\theta(-I)=1$. Write
$\rho_\theta=R_T^G(\theta)$ for its Deligne--Lusztig character.
Here $R_T^G(\theta)$ denotes the (a priori virtual) character
obtained by Deligne--Lusztig induction from $\theta$; in the
regular cases used below it is an actual irreducible character
because the Deligne--Lusztig sign is
$\varepsilon_G\varepsilon_T=(-1)^{2+0}=+1$, using the split ranks
of $G$ and $T$ \cite{DL}.

We use Srinivasan's labels~\cite{Srinivasan} for the unipotent
characters
\[
 \eta_+=\theta_9,\qquad\eta_-=\theta_{10},\qquad
 \gamma=\theta_{11},\qquad\xi=\theta_{12},\qquad
 \mathrm{St}=\theta_{13}.
\]
Their degrees, in this order, are
$P(P+1)^2/2$, $P(P-1)^2/2$, $P(P^2+1)/2$,
$P(P^2+1)/2$, and $P^4$.
For each orbit representative $\theta$ above,
$\rho_\theta$ is a regular cuspidal character of degree
$(P^2-1)^2$.

\begin{proposition}[Uniform irreducibility indicator]
\label{prop:generic-indicator}
As class functions on $G$ (and, equivalently, on $S$),
\begin{equation}\label{eq:indicator-general}
 m1_{\mathrm{irr}}
 =\kappa(1-\eta_++\eta_-+\mathrm{St})
  -\sum_{\theta\in\mathcal O}\rho_\theta.
\end{equation}
\end{proposition}
\begin{proof}
The Coxeter-torus Deligne--Lusztig decomposition is
\[
 R_T^G(1)=1-\eta_++\eta_-+\mathrm{St};
\]
see \cite[Table~2.8]{GeckMalle} and compare \cite{Srinivasan}.
Set
\[
 \Phi=\kappa R_T^G(1)
   -\sum_{\theta\in\mathcal O}R_T^G(\theta).
\]
We show that $\Phi(g)=m1_{\mathrm{irr}}(g)$ by examining the
semisimple part of $g\in G$.

If the semisimple part is not conjugate into $T$, the
Deligne--Lusztig character formula gives $\Phi(g)=0$.
Such an element cannot be irreducible on the natural module.

Suppose next that the semisimple part is central. Write $g=zu$
with $z\in\{\pm I\}$ and $u$ unipotent. The Jordan-decomposition
formula expresses $R_T^G(\theta)(zu)$ as the same Green-function
value multiplied by $\theta(z)$. As all characters used here are
trivial on the center, the $\kappa$ terms in the sum cancel the
first term. Hence $\Phi(zu)=0$; this includes $g=\pm I$.
No such $g$ acts irreducibly.

Finally, suppose the semisimple part is conjugate to an element
$t\in T\setminus\{\pm I\}$. Such an element is regular, so its
centralizer is $T$, and $g$ itself is conjugate to $t$.
The rational Weyl group has order four. Hence
\[
 R_T^G(1)(t)=4,
 \qquad
 \sum_{\theta\in\mathcal O}R_T^G(\theta)(t)
 =\sum_{\substack{\vartheta\in\widehat T_0\\\vartheta\ne1}}
   \vartheta(t_0)=-1,
\]
where $t_0$ is the image of $t$ in $T_0$.
The final equality is the orthogonality relation for characters
of the cyclic group $T_0$, evaluated at $t_0\ne1$.
It follows that
\[
 \Phi(t)=4\kappa+1=m.
\]
These, and only these, are the irreducible elements of $G$.
Thus $\Phi=m1_{\mathrm{irr}}$ on every conjugacy class.
\end{proof}

\begin{remark}[The density of irreducible elements]
\label{rem:irr-density}
As a consistency check, every irreducible element is conjugate
to one of the $P^2-1$ noncentral elements of $T$. Since each such
element has centralizer $T$ and rational Weyl group of order four,
\[
 \frac{|\{g\in G:1_{\mathrm{irr}}(g)=1\}|}{|G|}
 =\frac{P^2-1}{4(P^2+1)}=\frac{\kappa}{m}.
\]
This is also the coefficient of the trivial character in
\eqref{eq:indicator-general}.
\end{remark}

\begin{remark}\label{rem:three-indicator}
For $P=3$, one has $m=5$, $\kappa=1$, and $|\mathcal O|=1$.
The characters specialize to
\[
 \eta_-=\chi_6,\qquad \eta_+=\chi_{24},\qquad
 \rho_\theta=\chi_{64},\qquad \mathrm{St}=\chi_{81},
\]
where the subscripts indicate degrees. Consequently,
\begin{equation}\label{eq:indicator-three}
 5\,1_{\mathrm{irr}}=1+\chi_6-\chi_{24}-\chi_{64}+\chi_{81}.
\end{equation}
In particular, the character identity used in the explicit
enumeration for $P=3$ follows from the uniform Coxeter-torus
formula, rather than requiring an independent table calculation.
Summing it over the Frobenius classes of the parameters
$a\in\F_{3^H}$ expresses $5A_H(3)$ as four character sums.
Their cohomological evaluation below yields the explicit formula
in Theorem~\ref{thm:main}.
\end{remark}

Let $F_{a,H}$ denote the geometric Frobenius at the parameter
$a\in\F_{P^H}$ on the root local system. The fixed subspaces and the
irreducibility test are unchanged on passing from arithmetic Frobenius
to its inverse. We use geometric Frobenius throughout the trace formulas.
Define the formal power series (without assuming rationality at this stage)
\[
\mathcal R_P(u)=\exp\left(4m\sum_{H\geq1}A_H(P)\frac{u^H}{H}\right),
\qquad m=\frac{P^2+1}{2}.
\]

\begin{theorem}\label{thm:main}
Put
\begin{align*}
N(u)={}&1+u+2u^2+u^3+2u^4+17u^5+26u^6\\
 &+51u^7+18u^8+27u^9+162u^{10}+243u^{11}+729u^{12},\\
D(u)={}&1-u^2-3u^3+3u^4+18u^5+9u^6-27u^7-27u^8+243u^{10}.
\end{align*}
Then
\[
\mathcal R_3(u)=
 \left(\frac{N(u)}{(1-u)(1-3u)D(u)}\right)^4.
\]
Equivalently, if $N(u)=\prod_{i=1}^{12}(1-\alpha_i u)$ and
$D(u)=\prod_{j=1}^{10}(1-\beta_j u)$, then for every $H\geq1$,
\[
A_H(3)=\frac{3^H+1-\sum_i\alpha_i^H+\sum_j\beta_j^H}{5}.
\]
The cohomological input is established in
Sections~\ref{sec:monodromy} and~\ref{sec:three}. The exact
computation in Section~\ref{sec:certificate} then determines
the coefficients uniquely, using the exact finite-field counts stated below.
\end{theorem}

\subsection*{Notation}
The following symbols recur across the two parts.
\[
\begin{array}{ll}
 P=p^f,\ Q=P^H & \text{base and extension field sizes}\\
 G,\ S & \Sp_4(\F_P)\text{ and its central quotient}\\
 T,\ T_0 & \text{Coxeter torus and its image in }S\\
 B,\ W,\ U,\ U_0 & \text{root pairing, positive-valuation plane, and wild groups}\\
 \mathcal W & \text{total Weil representation space}\\
 \mathcal B_P,\ S_P(u),\ C_P(u) & \text{flag curve, Steinberg and cuspidal factors}\\
 N(u),\ D(u) & \text{the two factors at }P=3
\end{array}
\]

\section{Monodromy and ramification at infinity}\label{sec:monodromy}
\subsection{The Katz--Tiep construction and its relation to the root space}
We work for $P=p^f$, $p$ odd, $p\nmid f$, and fix an auxiliary prime
$\ell\ne p$. We briefly recall the particular part of the construction
of Katz--Tiep \cite{KT} needed here. Their local systems on the affine
parameter line are obtained from exponential sums and realize finite
Weil representations of symplectic groups. For $n=2$ the \emph{total}
Weil representation of $\Sp_4(\F_P)$ has dimension $P^2$ and is the
direct sum of its odd- and even-dimensional constituents, of dimensions
$(P^2+1)/2$ and $(P^2-1)/2$. The term ``total'' is important: using
both constituents retains the action of the central involution, which
will be needed to identify the kernel of the monodromy representation.

Here is the elementary bridge from the exponential sum to our skew
quartic. For a finite extension $k/\F_P$ and $b\in k$, consider the
quadratic form on the $\F_p$-vector space underlying $k$ given by
\[
 q_b(x)=\Tr_{k/\F_p}
 \bigl(x^{P^2+1}+x^{P+1}+b x^2\bigr).
\]
Its polar form has a radical that can be described by an additive
polynomial. More precisely, for each finite extension $k/\F_P$,
the radical consists of the elements of $k$ annihilated by this
polynomial: the polar equation, raised to the $P^2$-th power,
becomes the equation displayed below. In the normalization used below
this polynomial is
\[
 R_b(X)=X^{P^4}+X^{P^3}+2b^{P^2}X^{P^2}+X^P+X.
\]
Thus the purely inseparable change of parameter
\[
 a=2b^{P^2}
\]
identifies $R_b$ with $L_a$ from Section~\ref{sec:setup}. This is
the reason Katz--Tiep's sheaves are relevant to the present counting
problem: the quadratic form appearing in their trace function has
exactly the root space whose Frobenius action tests whether $f_a$
is an atom.

The second ingredient is the quadratic Gauss-sum identity. We record
its elementary content, to make the normalization of the subsequent
comparison transparent. Let $k/\F_P$ be finite, let
$\psi_k=\psi\circ\Tr_{k/\F_p}$, and set
\[
 S_b(k)=\sum_{x\in k}\psi_k
   (x^{P^2+1}+x^{P+1}+bx^2).
\]
Writing $q_b$ for the corresponding $\F_p$-valued quadratic form,
translation $(x,y)\mapsto(x+y,y)$ and additive-character
orthogonality give
\[
 |S_b(k)|^2
 =|k|\sum_{v\in\operatorname{rad}(q_b)}\psi(q_b(v))
 =|k|\,\#\operatorname{rad}(q_b).
\]
In the last equality, $q_b$ vanishes on its radical, since $p$ is
odd. The radical is exactly $\ker L_{2b^{P^2}}\cap k$.
Thus, for the normalization $A_k=A_{\F_p}^{[k:\F_p]}$ of \cite[\S3]{KT}, where
$A_{\F_p}=-\sum_{r\in\F_p^\times}\psi(2r)\chi_{2,\F_p}(r)$
because $P^2\equiv1\bmod4$, we have $|A_k|^2=|k|$ and
\begin{equation}\label{eq:gauss-kernel-direct}
 \left|-A_k^{-1}S_b(k)\right|^2
 =\#\{x\in k:L_{2b^{P^2}}(x)=0\}.
\end{equation}
This is the specialization of Katz--Tiep
\cite[Lemma~5.2]{KT} to $(n,s,t,q)=(2,1,b,P)$.
The arithmetic normalization of the total Weil trace function is
that of \cite[\S\S3,5]{KT}; the absolute-square identity does
not depend on the choice of phase of $A_k$. Since the trace of the
endomorphism representation of a finite unitary representation
is the squared absolute value of its character, the Frobenius
trace on the root permutation representation agrees with that
on the endomorphism of the total Weil sheaf. The latter is
insensitive to scalar operators, so the kernel comparison below is
needed to pass to \emph{full} symplectic monodromy.

\subsection{Full monodromy and the kernel comparison}
The relevant result is Katz--Tiep \cite[Theorem~3.2(iii)]{KT},
not a theorem about the natural four-dimensional representation.
For their specialization $s=1$, take $n=2$ and $q=P=p^f$.
For odd $p$, the theorem applies when $n\ge2$, $p\nmid n$,
and $q=p^f$ with $p\nmid f$. The parameter dictionary is
$(n,s,t,q)=(2,1,b,P)$: $s$ is fixed as the coefficient of
$x^{P+1}$, and $t=b$ varies as the coefficient of $x^2$ in the
total trace $-S_b(k)/A_k$. The theorem gives geometric monodromy
$\Sp(2n,q)$ on the direct sum of the two matched normalized
Fourier constituents, in its total Weil representation.
The system in the theorem starts on the parameter line over
$\F_p$; on its pullback to the line over $\F_q$, the arithmetic
and geometric groups agree. The two constituents must be kept together
in order to detect the central involution.

The change $a=2b^{P^2}$ is radicial and induces an equivalence
of the corresponding geometric finite \'etale sites. Under this
change the polar radical of the quadratic phase is the root space
of $L_a$, as made explicit in \cite[\S5, especially Lemma~5.2]{KT}.
It is important, however, that the monodromy theorem concerns the
Weil sheaf, whereas our counting problem concerns the natural
four-dimensional root-space representation. We therefore give the
additional kernel-comparison argument rather than identifying these
representations without proof.

\begin{lemma}[Kernel comparison in the monodromy argument]
\label{lem:monodromy-kernel}
Let $\Gamma$ be the geometric fundamental group of the parameter
line, let $\omega:\Gamma\twoheadrightarrow G$ be the Katz--Tiep
monodromy quotient, and let
\[
 \Omega:G\longrightarrow\operatorname{GL}(\mathcal W)
\]
be the total Weil representation, where $\dim \mathcal W=P^2$.
Let $\rho:\Gamma\to\Sp(V)$ denote the action on the four-dimensional
root space of $L_a$. Assume that, for every closed point and every
finite constant-field extension, the Frobenius traces of the associated
arithmetic representations satisfy
\[
 \operatorname{Tr}(\mathrm{Frob}\mid\overline{\mathbb Q}_\ell[V])
 =\operatorname{Tr}(\mathrm{Frob}\mid\operatorname{End}(\mathcal W)),
\]
where the latter is acted on by conjugation through $\Omega\circ\omega$.
Then $\rho$ factors through $G$, and the induced map
$G\to\Sp(V)$ is an isomorphism.
\end{lemma}
\begin{proof}
The permutation representation on the finite set $V$ has Frobenius
trace equal to the number of fixed roots. The action of $g\in G$ on
$\operatorname{End}(\mathcal W)$ is
\[
 T\longmapsto\Omega(g)T\Omega(g)^{-1}.
\]
Via $\operatorname{End}(\mathcal W)\simeq \mathcal W\otimes \mathcal W^\vee$, its trace is
\[
 \operatorname{Tr}(\Omega(g))\operatorname{Tr}(\Omega(g)^{-1})
 =|\operatorname{Tr}(\Omega(g))|^2,
\]
where the final equality uses a complex realization of the finite
Weil representation. The quadratic Gauss-sum identity
\cite[Lemma~5.2]{KT}, after $a=2b^{P^2}$, identifies this with
the number of Frobenius-fixed roots. The finite-image normalization and arithmetic extensions are those
in \cite[\S\S3,5]{KT}: one must compare traces of the \emph{same}
normalized total Weil sheaf in every finite extension. This avoids
inferring equality of arithmetic representations from a comparison
of geometric traces alone.

More formally, apply the comparison first to the arithmetic
fundamental group. Chebotarev density gives equality of characters
from equality of closed-point Frobenius traces. Both representations
have finite image and are semisimple in characteristic zero; thus
Brauer--Nesbitt makes them isomorphic. Restricting to $\Gamma$ gives
\[
 \overline{\mathbb Q}_\ell[V]\simeq
 \operatorname{End}(\mathcal W)
 \qquad\text{as $\Gamma$-representations}.
\]
The kernel of the first representation is $\ker\rho$, since the
permutation action on the basis indexed by $V$ is faithful. The kernel
of the second is the inverse image, under $\omega$, of the projective
kernel of $\Omega$, because conjugation on all endomorphisms is trivial
exactly for scalar operators.

The total Weil representation is projectively faithful. By its
\emph{projective kernel} we mean
\[
 \ker(G\longrightarrow\operatorname{PGL}(\mathcal W))
 =\{g\in G:\Omega(g)\text{ acts as a scalar on }\mathcal W\}.
\]
It is precisely the kernel of the conjugation action on
$\operatorname{End}(\mathcal W)$, and it is a normal subgroup of $G$. The central
quotient $\operatorname{PSp}_4(\F_P)$ is simple for odd $P$;
together with the perfectness of $G$, this shows that a proper normal
subgroup of $G$ must be central. The kernel is proper since the
projective Weil representation is nontrivial. Moreover $-I$ acts by
opposite scalars on the two nonzero even and odd Weil constituents,
so it is not scalar on their direct sum $\mathcal W$. Thus the projective
kernel is trivial. Consequently
\[
 \ker\rho=\ker\omega.
\]
Thus $\rho$ and $\omega$ have the same geometric kernel, and
$\rho$ descends through the quotient $\Gamma/\ker\omega\simeq G$.
The resulting homomorphism $G\hookrightarrow\Sp(V)$ is injective.
Both groups have order $P^4(P^2-1)(P^4-1)$, so it is an isomorphism.
In particular, the natural root representation and the total Weil
representation have the same monodromy quotient, up to this
identification.
\end{proof}

\subsection{The small Weil constituent and local Fourier transformation}
To compute the local ramification it is enough to work with one
\emph{faithful} constituent of the total Weil representation. We spell
out why this is legitimate and why the binomial Artin--Schreier sheaf
appears.

Set
\[
 d=\frac{P^2+1}{2},\qquad e=\frac{P+1}{2},\qquad
 k=d-1=\frac{P^2-1}{2},\qquad F(x)=x^d+x^e.
\]
The total Weil representation decomposes into constituents of dimensions
$d$ and $k$. In the Fourier--Deligne realization of Katz--Tiep
\cite{KT}, the rank-$k$ (small) constituent is, up to a constant-field
rank-one normalization, the Fourier transform of the rank-one
Artin--Schreier sheaf
\[
 \mathcal L_\psi(F(x)).
\]
The complementary constituent is obtained by inserting the quadratic
Kummer character in the input variable before Fourier transformation.
At the level of trace functions, for a finite extension $K/\F_P$,
these two transforms involve respectively the sums
\[
 \sum_{x\in K}\psi_K(F(x)+bx),\qquad
 \sum_{x\in K^\times}\chi_{2,K}(x)\psi_K(F(x)+bx).
\]
Indeed the elementary square-map identity
\[
 \sum_{y\in K}h(y^2)
 =\sum_{x\in K}(1+\chi_{2,K}(x))h(x),
 \qquad \chi_{2,K}(0)=0,
\]
explains why their sum is the total Weil trace, whose phase is
\[
 F(y^2)+by^2=y^{P^2+1}+y^{P+1}+by^2.
\]
Thus the binomial $F$ is not an auxiliary polynomial introduced for
the ramification calculation: it is the defining phase of the
small constituent in Katz--Tiep's construction.

The small Weil representation $\Omega_{\rm sm}:G\to
\operatorname{GL}(\mathcal W_{\rm sm})$ is faithful. One can see this from the
normal-subgroup structure of $G$: its kernel is proper and hence
central, whereas the central involution acts nontrivially on the
small constituent.

We use the upper numbering of ramification groups in the convention
of Serre. For $v\ge0$, let $I_\infty^v$ denote the upper-numbering
subgroup of geometric inertia at infinity; for a finite-image
representation, its positive upper breaks are the numbers $b>0$ at
which the image filtration jumps, or equivalently where its image at
$b$ differs from the image at $b+$ (the subgroup for $v>b$
sufficiently close to $b$). The slopes of an individual representation
need not detect every break of the image filtration. Write
\[
 J^v=\omega(I_\infty^v)\subset G,\qquad
 J^{b+}=\bigcup_{\varepsilon>0}J^{b+\varepsilon}.
\]
The filtration is decreasing; the union defines the right-hand
subgroup because its image is finite. Faithfulness gives
\[
 J^u\xrightarrow{\ \sim\ }
 \Omega_{\rm sm}(J^u)
 \qquad (u\ge0),
\]
compatibly with subgroup inclusions. By the preceding kernel
comparison the same $J^u$ acts faithfully on the natural root space.
Hence the ranks, orders, and jumps of the \emph{images} of upper
ramification groups may be calculated on the small Weil sheaf.
This does not by itself locate their images inside a specified
Siegel parabolic; that identification will be made using the roots.

The next calculation is the main local input to the paper. It is
useful to separate three steps: first compute the polar phases after
a tame pullback, then recover the two filtered pieces of the abelian
wild image from their characters, and finally identify that image
with the Siegel radical by an independent Newton-polygon argument.
The last step is necessary because local Fourier analysis alone gives
an abstract filtered group, not its embedding in $\Sp_4(\F_P)$.

We use Laumon's local Fourier transform and stationary-phase
formalism \cite{Laumon} and the explicit local Fourier formula of
Abbes--Saito \cite[Theorem~3.9]{AS}. The latter, under its Legendre
conditions, describes the transformed local representation as the
pushforward of a rank-one Artin--Schreier character (with tame and
unramified corrections). It is this explicit pushforward description,
rather than the mere existence of local Fourier transformation, that
allows us to read the ramification filtration.

\begin{remark}[Independent check from Laumon's local Fourier theory]
\label{rem:laumon-rank-check}
Here is a check of the normalization and the largest ramification
break that does not depend on the explicit stationary-phase expansion
below.  Laumon's local Fourier functor
$\mathcal F^{(\infty,\infty')}$ takes a local representation of
slope greater than one to one of slope greater than one; by
\cite[Th\'eor\`eme~2.4.3(iii)(b)]{Laumon}, its rank and Swan
conductor are, respectively,
\[
 \operatorname{rank}(\mathcal F^{(\infty,\infty')}(V))
   =\operatorname{Sw}(V)-\operatorname{rank}(V),\qquad
 \operatorname{Sw}(\mathcal F^{(\infty,\infty')}(V))
   =\operatorname{Sw}(V).
\]
For $V=\mathcal L_\psi(x^d+x^e)$ near infinity, one has
$\operatorname{rank}(V)=1$ and $\operatorname{Sw}(V)=d$,
since $d>e$ and $p\nmid d$.  Consequently its local Fourier
transform has rank $d-1=k$ and a single slope $d/k$.  In
particular, this predicts the rank of the small Weil constituent
and its largest upper break
$\frac{d}{k}=\frac{P^2+1}{P^2-1}$.
Laumon's discussion of stationary-phase expansions
\cite[\S2.6.3]{Laumon} also explains the Legendre relation
$c=-F'(x)$; in the rank-one case he records the additional
assumptions $p\nmid d(d-1)$, both satisfied here.
The assertion that the small Weil representation has a single slope
does not imply that its finite wild inertia image has only one
upper-numbering jump: the former describes the breaks of the given
representation, whereas the latter describes the complete filtration
of the group image. This independent check does \emph{not}
determine the smaller break $e/k$ or the intermediate wild subgroup
$U_0$: these require the explicit Abbes--Saito formula and the
character-space calculation in the proof below.
\end{remark}

Recall that a \emph{Singer subgroup} of
$\operatorname{GL}_2(\F_P)$ is the image of
$\F_{P^2}^{\times}$ acting by multiplication on the
$2$-dimensional $\F_P$-space $\F_{P^2}$.
It is cyclic of order $P^2-1$, and a generator is called a
\emph{Singer cycle}. We will use this terminology for its
conjugates in $\operatorname{GL}(W)$ as well.

\begin{lemma}[Direct calculation of tame inertia on the roots]
\label{lem:tame-roots}
Let $K=\overline{\F}_P((u))$, where $u=a^{-1}$, and let $V$
be the four-dimensional root space of $L_a$. The roots of strictly
positive valuation, together with zero, form an inertia-stable
$2$-dimensional $\F_P$-subspace $W$. Wild inertia fixes $W$
pointwise, and the inertia image on $W$ is a Singer subgroup
of order exactly $P^2-1$. The plane $W$ is Lagrangian, and the
prime-to-$p$ quotient of the full inertia image on $V$ has order
$P^2-1$.
\end{lemma}
\begin{proof}
With $v(u)=1$, the polynomial whose nonzero roots we study is
\[
 \frac{L_a(X)}X
 =X^{P^4-1}+X^{P^3-1}+u^{-1}X^{P^2-1}+X^{P-1}+1.
\]
Its lower Newton polygon has two segments, joining $(0,0)$ to
$(P^2-1,-1)$ and then to $(P^4-1,0)$. Their respective slopes
are $-1/(P^2-1)$ and $1/(P^4-P^2)$. Thus precisely
$P^2-1$ nonzero roots have valuation $1/(P^2-1)$; the remaining
$P^4-P^2$ nonzero roots have valuation $-1/(P^4-P^2)$. Because $L_a$ is additive,
the zero root together with the roots of positive valuation is
closed under addition and multiplication by $\F_P$. It therefore
forms an inertia-stable subspace $W\subset V$ of cardinality $P^2$.

Choose $\lambda$ with $\lambda^{P^2-1}=u$ and substitute
$X=\lambda Y$. After division by $\lambda$ we obtain
\[
 Y+Y^{P^2}+\lambda^{P-1}Y^P+
 \lambda^{P^3-1}Y^{P^3}+\lambda^{P^4-1}Y^{P^4}=0.
\]
The reduction modulo $\lambda$ is $Y^{P^2}+Y=0$. Its derivative
is $1$, so its $P^2$ roots in the algebraically closed residue
field are simple. Hensel's lemma lifts each root uniquely to
$\overline{\F}_P((\lambda))$. Every element of $W$ is among these
lifts, because its rescaled coordinate $Y=X/\lambda$ has
nonnegative valuation. Thus $W$ is pointwise fixed by wild
inertia: the extension $K(\lambda)/K$ is tame of degree $P^2-1$.

The residual root space
\[
 W_{\mathrm{red}}=\{y\in\overline{\F}_P:y^{P^2}+y=0\}
\]
is one-dimensional over $\F_{P^2}$. A tame inertia element taking
$\lambda$ to $\zeta\lambda$, with $\zeta\in\mu_{P^2-1}$,
acts on this residual root space by multiplication by $\zeta$.
Indeed, if $g(\lambda)=\zeta\lambda$, then
\[
 \overline{g(X)/\lambda}
 =\overline{(g(\lambda)/\lambda)\,g(X/\lambda)}
 =\zeta\,\overline{X/\lambda},
\]
because geometric inertia acts trivially on the residue field.
These multiplications give the entire $\F_{P^2}^\times$-action,
and hence an injective cyclic subgroup of
$\operatorname{GL}_{\F_P}(W)$ of order $P^2-1$ (a Singer subgroup).
The inertia action on $W$ factors through this group because
all elements of $W$ are defined over $K(\lambda)$.

The alternating form on $V$ is inertia invariant. If its
restriction to the two-dimensional space $W$ were nonzero, it
would be nondegenerate, so the inertia image on $W$ would lie in
$\operatorname{Sp}(W)=\operatorname{SL}_2(\F_P)$. But the
determinant of multiplication by $\zeta\in\F_{P^2}^\times$
is its norm to $\F_P^\times$, and this norm is surjective.
Hence the restriction of the alternating form must vanish.
As $\dim W=\frac12\dim V$, the plane $W$ is Lagrangian.

Consequently the full inertia image $I$ lies in the Siegel
parabolic $P_W=\operatorname{Stab}_{\Sp(V)}(W)$. The restriction map
$P_W\to\operatorname{GL}(W)$ has kernel the Siegel unipotent radical
$U_{\mathrm{Sieg}}$. The kernel of its restriction to $I$ is
$I\cap U_{\mathrm{Sieg}}$, a $p$-group. The image
of $I$ in $\operatorname{GL}(W)$ is the Singer group of order
$P^2-1$, while the wild inertia image is contained in
$U_{\mathrm{Sieg}}$. Since the kernel of
$I\to\operatorname{GL}(W)$ is a $p$-group and the quotient has
order prime to $p$, this kernel is exactly the wild inertia
image (the normal Sylow $p$-subgroup of $I$). Thus the tame
quotient has order $P^2-1$, as claimed.
\end{proof}

For subsequent use, let
\[
 \operatorname{Sym}_2(\F_P)
 =\{B\in M_2(\F_P):B^{\mathsf t}=B\}.
\]
This is a three-dimensional $\F_P$-vector space, considered
as an additive group. In a symplectic basis adapted to a
Lagrangian plane, the Siegel unipotent radical is
\[
 U_{\mathrm{Sieg}}=
 \left\{
 \begin{pmatrix}I_2&B\\0&I_2\end{pmatrix}:
 B\in\operatorname{Sym}_2(\F_P)
 \right\}.
\]
The multiplication law is $u_Bu_C=u_{B+C}$; thus this is an
isomorphism of additive groups, not an assertion about matrix
multiplication in $\operatorname{Sym}_2(\F_P)$.

\begin{lemma}[Wild inertia and its breaks]\label{lem:wild-image}
The wild inertia image $U$ on the natural root space is elementary
abelian of order $P^3$. Its upper breaks are $b_1=1/(P-1)$ and
$b_2=(P^2+1)/(P^2-1)$. Between these breaks its image is a subgroup
$U_0$ of order $P^2$.
\end{lemma}
\begin{proof}
By the preceding discussion, the rank-$k$ small Weil sheaf is the
Fourier--Deligne transform of $\mathcal L_\psi(F(x))$ up to an
unramified constant-field normalization. Its trace over $K$ is
\[
 -A_K^{-1}\sum_{x\in K}\psi_K(F(x)+bx),
 \qquad A_K=A_{\F_p}^{[K:\F_p]}.
\]
The normalization does not change geometric inertia. The parameter change $a=2b^{P^2}$ is a scalar change followed by
an iterate of relative Frobenius. On the completed local fields at
infinity, it identifies the ramification filtrations: the scalar change
preserves the normalized valuation, and the purely inseparable part
induces the same map on finite separable extensions as the Frobenius
identification, which preserves normalized valuations with the
corresponding uniformizers. The local Fourier calculation therefore
applies to the inertia quotient in Lemma~\ref{lem:monodromy-kernel}.

Apply Abbes--Saito \cite[Theorem~3.9]{AS} at infinity to the
rank-one sheaf $\mathcal G=\mathcal L_\psi(F(x))$. Their hypotheses
in this case are (i) the Legendre triple condition, (ii) finite
source and target maps that are \'etale on generic points, (iii)
a source point at infinity, and (iv$'$) a target point at infinity
and all slopes of the source pushforward strictly greater than one.
We work over the henselization at infinity with perfect residue
field $\overline{\F}_P$, in characteristic $p>2$, and use its
completion for the valuation calculation.

In the uniformizer $v=1/x$, set
\[
 A=v^{-d}+v^{-e},\quad b=v^{-1},\quad
 c=-\tfrac12(v^{-k}+v^{-(e-1)}),\qquad
 \alpha=A'=-d v^{-d-1}-e v^{-e-1},
\]
with primes denoting $d/dv$. Here $d=e=1/2$ and $k=-1/2$ as
elements of the residue field. In Definition~2.16 of \cite{AS},
take the tame factor to be trivial and the wild factor to be
the Artin--Schreier sheaf of $A$. It is trivialized by a cyclic
extension of degree $p$, since $A$ has pole order $d$ prime to $p$.
For this length-one Witt vector, the conditions (2.13.1)--(2.13.3)
of \cite{AS} are
\[
 A\in\mathrm{fil}_n,\qquad
 2\operatorname{ord}_v(\alpha+cb')\ge -n+\nu(c),\qquad
 2\nu(b)+p\nu(c)<(p-2)n,
\]
where $n=-\operatorname{ord}_v(v\alpha)$ and
$\nu(h)=\operatorname{ord}_v(vh'/h)$. The pole orders give $n=d$
and $A\in\mathrm{fil}_d$. We have $\alpha+cb'=0$ and
$\nu(b)=\nu(c)=0$, since the leading exponents $-1$ and $-k$ are
prime to $p$. Hence the last inequality is $0<(p-2)d$.
The source map $b=x$ has degree one and the target map $c=-F'(x)$
has local degree $k$ prime to $p$. Both are finite on the
henselian traits, generically \'etale, and have poles at $v=0$.
The source pushforward has the single slope $d>1$.
The formula therefore gives the pushforward along $c$ of the
character with phase $A+bc=F-xF'$ up to tame factors.

Choose a tame coordinate $z$ with $c=-z^k/2$.
Writing $h=d-e=P(P-1)/2$, we have
 $z^k=x^k(1+x^{-h})$, hence
 $z=x(1+x^{-h})^{1/k}$ and
 $x=z+2z^{1-h}+O(z^{1-2h})$; here $1/k=-2$ in
 characteristic $p$.
 Since $d-2h<0$ and $e-h<0$ for $P\ge3$,
 \[
 x^d=z^d+z^e+O(z^{d-2h}),\qquad
 x^e=z^e+O(z^{e-h}).
 \]
 As $1-d=1-e=1/2$ in characteristic $p$,
 $F-xF'=\tfrac12(x^d+x^e)$, and therefore its polar part is exactly
\[
F(x)-xF'(x)=\tfrac12z^d+z^e\pmod{z^{-1}\overline{\F}_P[[z^{-1}]]}.
\]
After the tame degree-$k$ pullback, its conjugate wild characters are
\[\mathcal L_\psi\bigl(\tfrac12\zeta z^d+\zeta^e z^e\bigr),
\qquad \zeta\in\mu_k.\]
We now prove the linear independence assertion used to determine the
wild image, rather than only citing the two coefficient fields.
Fix a primitive $k$th root of unity $\zeta$ and recall that
$k=e(P-1)$ and $d=k+1$.  The conjugate Artin--Schreier
classes have pole terms
\[
 \tfrac12\zeta^j z^d+\zeta^{ej}z^e,
 \qquad j\in\mathbb Z/k\mathbb Z.
\]
Because $p\nmid de$ and $d\ne e$, a nonzero linear combination of
$z^d$ and $z^e$ cannot be an Artin--Schreier coboundary: the
pole order of a nontrivial coboundary with a pole is divisible by $p$.
Thus their $\F_p$-linear relations are exactly the relations between
the coefficient pairs
\[
 (\zeta^j,\zeta^{ej})\in\F_{P^2}\oplus\F_P.
\]
The scalar $1/2$ in the first component is invertible and hence
irrelevant to their span.  The order of $p$ modulo $k$ is $2f$:
it divides $2f$ because $k\mid(P^2-1)$, and it cannot be a proper
divisor of $2f$, since every proper divisor is at most $f$ and
$k>P=p^f$.  Similarly, the order of $p$ modulo $P-1$ is $f$.
Consequently, the minimal polynomials $m_d(T)$ and $m_e(T)$ of
$\zeta$ and $\zeta^e$ over $\F_p$ are distinct irreducibles of
degrees $2f$ and $f$.  By the Chinese remainder theorem, evaluation
induces a surjection
\[
 \F_p[T]\longrightarrow\F_{P^2}\oplus\F_P,
 \qquad h(T)\longmapsto(h(\zeta),h(\zeta^e)),
\]
whose kernel is $(m_dm_e)$.  Since $\zeta^k=1$, every evaluation
is in the span of the $k$ pairs displayed above.  Their span
therefore has dimension $2f+f=3f$, and, as a module for the tame
cyclic quotient, splits into the two distinct irreducible summands
coming from $m_d$ and $m_e$.

After the tame degree-$k$ pullback, the wild inertia representation
is the direct sum of the above rank-one characters.  Let $U$ be its
finite wild image and let $M\subseteq U^\vee$ be the subgroup
generated by these characters.  They separate points of $U$ because
$U$ is defined as the image on the small Weil constituent; hence
$M=U^\vee$.  The calculation just made gives
\[
 M=M_d\oplus M_e,\qquad
 \dim_{\F_p}M_d=2f,\quad \dim_{\F_p}M_e=f.
\]
In particular, $U^\vee$ is elementary abelian of order $p^{3f}$,
so $U$ is elementary abelian of order $P^3$.

The conductor of a nontrivial class in $M_e$ is $e$, while a class
with nonzero $M_d$-component has conductor $d$: the nonzero leading
pole cannot be removed by an Artin--Schreier coboundary, since
$p\nmid de$. Thus the lower-break characters form exactly $M_e$.

For a finite abelian $p$-group $U$ and any subgroup $H\subset U$,
Pontryagin duality recovers $H$ as the intersection of the kernels
of all characters that are trivial on $H$. The subgroup $M_e$ consists exactly of the characters whose
restriction has conductor at most the smaller of the two breaks.
Its common kernel therefore detects the part of $U$ still visible
after the first break. Applying this to the upper ramification
groups, the subgroup surviving between the two breaks is
\[
 U_0=\bigcap_{\chi\in M_e}\ker\chi,
 \qquad (U/U_0)^\vee=M_e.
\]
In particular
\[
 |U|=|M|=p^{3f}=P^3,\qquad
 |U/U_0|=|M_e|=p^f=P,\qquad |U_0|=P^2.
\]
One should distinguish this from the fact that each individual
rank-one summand can have conductor $d$: the smaller break is detected
by products of diagonal characters in which the highest pole cancels.
Finally, upper breaks are divided by $k$ when descending from the
tame degree-$k$ extension. Thus the two positive breaks on the
original parameter line are
\[
 \frac ek=\frac1{P-1},\qquad
 \frac dk=\frac{P^2+1}{P^2-1},
\]
and the successive positive-level images are $U$, $U_0$, and $1$.
 The faithfulness argument preceding this lemma now transfers the
filtration of the image of wild inertia to the natural root space.
In particular, it transfers the intermediate subgroup and both breaks,
not just the order of the full wild image.

\end{proof}

\begin{lemma}[Inertia in the Siegel parabolic]\label{lem:inertia}
The inertia image is $I=U\rtimes\langle s\rangle$, where
$U=U_{\mathrm{Sieg}}$ and $s$ acts on $W$ as a Singer cycle
of order $P^2-1$. Under the Singer action on $U$, the subgroup
$U_0$ of Lemma~\ref{lem:wild-image} is the two-dimensional
summand of $\operatorname{Sym}_2(\F_P)$.
\end{lemma}
\begin{proof}
The Fourier-transform computation above, together with the
faithfulness and kernel comparison, shows that the image $U$
of wild inertia on the natural root space has order $P^3$.
Lemma~\ref{lem:tame-roots} places it inside the Siegel radical
$U_{\mathrm{Sieg}}$, which also has order $P^3$ because
\[
 U_{\mathrm{Sieg}}\simeq
 \operatorname{Sym}_2(\F_P)
\]
as an additive group. Therefore $U=U_{\mathrm{Sieg}}$.
The tame Singer group has order prime to $p$; the resulting
extension splits by Schur--Zassenhaus. In particular, the full
inertia image is, up to conjugacy in $G$,
\[
 I=U\rtimes\langle s\rangle,
 \qquad |U|=P^3,\quad |s|=P^2-1.
\]
This identifies the inertia subgroup as a concrete subgroup of
$\Sp_4(\F_P)$, rather than only determining its abstract order.

The character calculation in Lemma~\ref{lem:wild-image} gives
the tame-stable filtration of the wild image
\[
 U\supset U_0\supset\{1\},
 \qquad |U_0|=P^2,
\]
where $U_0$ is the annihilator of the lower-break character
space $M_e\subset U^\vee$. In particular, $U_0$ is not
obtained from the Newton polygon alone: its order and break
come from the local Fourier computation. Under the Singer
action on $U\simeq\operatorname{Sym}_2(\F_P)$, the
three-dimensional $\F_P$-space decomposes into a
one-dimensional norm-character summand and a two-dimensional
summand, the latter being irreducible under the Singer cycle.
These summands are irreducible and nonisomorphic even over $\F_p$.
Indeed, if the Singer generator has eigenvalue $\beta$ on $W$ over
$\F_{P^2}$, the eigenvalues on symmetric forms are $\beta^{P+1}$
for the norm summand and $\beta^2,\beta^{2P}$ for the other
summand. Their respective minimal-polynomial degrees over $\F_p$
are $f$ and $2f$ by the same order calculation used above for
$P-1$ and $k$. As the Singer group has order prime to $p$, its
only invariant subgroup of order $P^2$ in $U$ is the latter
summand, which therefore equals $U_0$. This explains the identification
of $U_0$ used in subsequent character averages.

\end{proof}

For a subgroup $J$ and a representation $E$ of $G$, write $E^J$
for its fixed subspace and $\Sw_\infty(E)$ for the Swan conductor
at infinity.

\begin{lemma}[Swan conductor]\label{lem:swan}
For every characteristic-zero representation $E$ of $G$,
\begin{equation}\label{eq:swan}
 \Sw_\infty(E)=
 \frac{\dim E-\dim E^U}{P-1}
 +\frac{P}{P+1}(\dim E-\dim E^{U_0}).
\end{equation}
\end{lemma}
\begin{proof}
By Lemma~\ref{lem:wild-image}, the successive positive-level
images are $U$, $U_0$, and $1$. Integrating
$\dim E-\dim E^{J^v}$ over positive upper levels gives
$b_1(\dim E-\dim E^U)+(b_2-b_1)(\dim E-\dim E^{U_0})$.
Here $b_1=1/(P-1)$ and $b_2-b_1=P/(P+1)$.
\end{proof}

\subsection{Cohomological factors}
The purpose of the preceding inertia computation is to determine
the cohomological degrees of the trace factors used below.
Every representation of $S=G/\{\pm I\}$, pulled back
along the monodromy quotient, gives an $\ell$-adic local system
on the parameter line. For a character $\chi$ of $S$, denote it
by $\mathcal E_\chi$, with coefficients in
$\overline{\mathbb Q}_\ell$ for $\ell\ne p$, and set
\[
 L_\chi(u)=\det\bigl(1-uF\mid
 H^1_c(\mathbb A^1_{\overline{\F}_P},\mathcal E_\chi)\bigr),
\]
where $F$ is the geometric Frobenius of the base field.
For a nontrivial irreducible $\chi$, full geometric monodromy
implies $H_c^0=H_c^2=0$. The Grothendieck trace formula then
gives
\[
 \sum_{a\in\F_{P^H}}\chi(F_{a,H})
 =-\Tr\bigl(F^H\mid
 H^1_c(\mathbb A^1_{\overline{\F}_P},\mathcal E_\chi)\bigr).
\]
Thus the finite-field character sums are encoded by the
Frobenius polynomial $L_\chi$. The Swan conductor formula
\eqref{eq:swan}, together with Grothendieck--Ogg--Shafarevich,
will determine the degrees of these polynomials.

\section{Cohomological degrees and duality for \texorpdfstring{$P=3$}{P=3}}\label{sec:three}
In this section $P=3$, so $S\simeq \operatorname{PSU}_4(2)$. Denote by $\chi_6,\chi_{24},
\chi_{64},\chi_{81}$ its irreducible characters with the indicated
degrees. The last two are the Steinberg characters in the unitary and
symplectic descriptions respectively. By the specialization of the uniform Deligne--Lusztig identity
\eqref{eq:indicator-general} (see Remark~\ref{rem:three-indicator}),
\[
 5\,1_{\rm irr}=1+\chi_6-\chi_{24}-\chi_{64}+\chi_{81}.
\]
In $S$, irreducibility is equivalent to membership in its unique
order-five class.

The purpose of this identity is to convert the condition that a
Frobenius element act irreducibly into character sums. Indeed,
\[
 A_H(3)=\sum_{a\in\F_{3^H}}1_{\mathrm{irr}}(F_{a,H}),
\]
so the right-hand side of \eqref{eq:indicator-three} can be treated
by the Grothendieck trace formula. We first compute the dimensions
of the inertia-invariant subspaces needed for this passage.

Identify $U$ with symmetric $2\times2$ matrices $B$ over $\F_3$ via
$u_B=\left(\begin{smallmatrix}1&B\\0&1\end{smallmatrix}\right)$.
Choose the Singer cycle on $W$ as
$A=\left(\begin{smallmatrix}1&2\\1&1\end{smallmatrix}\right)$;
$A^{\mathsf t}A=2I$ and $U_0$ is the trace-zero plane.
The needed values and numbers of elements are
\[
\begin{array}{c|rrrr}
& B=0&\operatorname{rank}B=1&\det B=1&\det B=2\\\hline
\#\text{ in }U&1&8&6&12\\
\#\text{ in }U_0&1&0&4&4\\
\chi_6&6&-3&3&0\\
\chi_{24}&24&6&0&3\\
\chi_{64}&64&-8&4&-2\\
\chi_{81}&81&0&0&0
\end{array}
\]
The rank-two class labels can be checked without a naming convention:
the trace on the small Weil representation is $(S_B-1)/2$, where $S_B$
is the two-variable quadratic Gauss sum. It is $-2$ for $\det B=1$ and
$1$ for $\det B=2$. The character of its exterior square, $\chi_6$,
is therefore $3$ and $0$ respectively; for rank one it is $-3$.

The rank-one entries aggregate two congruence classes: for
$2\times2$ symmetric matrices of odd characteristic, congruence
classes are distinguished by rank and, in each positive rank, by
the discriminant square class. The two rank-one classes give the
same displayed character values, so they have been combined in
one column. The two nonsingular columns are distinguished by
the determinant square class in $\F_3^{\times}$.

For a representation $E$ of a finite group $J$, averaging the
operators in $J$ gives the projector onto $E^J$. Hence
\[
 \dim E^J=\frac{1}{|J|}\sum_{g\in J}\chi_E(g).
\]
For example, the first two rows of the table yield
\[
 \dim E_6^U=\frac{6+8(-3)+6(3)+12(0)}{27}=0,
 \qquad
 \dim E_6^{U_0}=\frac{6+4(3)+4(0)}9=2.
\]
The remaining invariant dimensions follow by the same averaging
calculation. We then insert them into the Swan conductor formula
\eqref{eq:swan}.

We write
$\mathcal E_j$ and $L_j$ for the objects attached to $\chi_j$.
For these nontrivial characters, full geometric monodromy
implies $H^0_c=H^2_c=0$. Since the local systems are lisse on the
affine line, Grothendieck--Ogg--Shafarevich gives
\[
 \chi_c(\mathbb A^1,\mathcal E_j)
 =\dim E_j-\Sw_\infty(E_j).
\]
Consequently,
\[
 \dim H^1_c(\mathbb A^1,\mathcal E_j)
 =\Sw_\infty(E_j)-\dim E_j.
\]
For instance, when $j=64$ the table above gives
$\dim E_{64}^U=0$ and $\dim E_{64}^{U_0}=8$, whence
\[
 \Sw_\infty(E_{64})
 =\frac{64}{2}+\frac34(64-8)=74,
 \qquad \dim H^1_c(\mathbb A^1,\mathcal E_{64})=10.
\]
The same calculation for the other three characters yields
\begin{equation}\label{eq:dimensions-three}
\begin{array}{c|rrrr}
j&\dim E_j^U&\dim E_j^{U_0}&\Sw_\infty(E_j)&
\dim H^1_c(\mathbb A^1,\mathcal E_j)\\\hline
6&0&2&6&0\\24&4&4&25&1\\64&0&8&74&10\\81&3&9&93&12
\end{array}
\end{equation}
The last column determines the degrees of the cohomological
Frobenius polynomials; it does not yet determine their coefficients.
Thus $L_6=1$. At $H=1$, $a=0,1,2$ give projective
orders $9,5,12$, and the $\chi_{24}$ values are $0,-1,0$.
The trace formula gives $L_{24}=1-u$.

We next justify the reciprocal symmetries used to determine $L_{64}$
and $L_{81}$. These symmetries reduce the number of coefficients that
must be recovered from the exact finite-field counts in the next
section. For $E_{64}$, \eqref{eq:dimensions-three} already implies $E_{64}^I=0$.
For $E_{81}$, the $U$-invariants of the symplectic Steinberg representation
are the Steinberg representation of the Levi $\mathrm{GL}_2(3)$.
The Singer subgroup of order eight has trace $3$ at its two scalar
elements and $-1$ at its other six elements. Its invariant dimension
is $(2\cdot3-6)/8=0$. Equivalently, this follows by averaging the
Steinberg character on the four cosets in $I/\{\pm1\}U$.
Hence both local systems have no inertia invariants at infinity.
Their compactly supported first cohomology equals their middle-extension
cohomology and is pure of weight one. The Steinberg representation
has a symmetric invariant form, for example from its real realization
in the building. The regular cuspidal character $\chi_{64}$ is
self-dual: the contragredient of $R_T^G(\theta)$ has parameter
$\theta^{-1}$, and $\theta^{-1}=\theta^{P^2}$ is in the same
rational Weyl orbit. Poincar\'e duality therefore pairs the Frobenius
eigenvalues of each cohomology group by $\alpha\mapsto3/\alpha$.
For the Steinberg factor, the induced cohomology pairing is alternating
with Frobenius multiplier three, so the determinant has positive sign.
For the cuspidal factor we retain both signs initially. Writing
$L_{81}=\sum n_i u^i$ and $L_{64}=\sum d_i u^i$ gives
\begin{equation}\label{eq:reciprocity-three}
 n_{12-i}=3^{6-i}n_i\ (0\leq i\leq6),\qquad
 d_{10-i}=\varepsilon\,3^{5-i}d_i\ (0\leq i\leq5),
 \qquad \varepsilon\in\{1,-1\}.
\end{equation}
The finite reconstruction in Section~\ref{sec:certificate}
will also determine $\varepsilon$.

Finally, sum \eqref{eq:indicator-three} over all parameters
$a\in\F_{3^H}$. The trivial character contributes $3^H$,
while $L_6=1$ makes the $\chi_6$-sum vanish. Since
$L_{24}=1-u$, the $\chi_{24}$-sum equals $-1$.
For $j=64,81$, the trace formula identifies the corresponding
character sum with the negative of the Frobenius trace on $H_c^1$.
We obtain
\begin{equation}\label{eq:trace-three}
5A_H(3)=3^H+1-\Tr(F^H\mid H^1_c(\mathbb{A}^1,\mathcal E_{81}))
+\Tr(F^H\mid H^1_c(\mathbb{A}^1,\mathcal E_{64})).
\end{equation}
In particular the degrees $12$ and $10$ in this calculation are theoretical
dimensions, not degrees selected by a fitting procedure.

\section{Exact reconstruction for \texorpdfstring{$P=3$}{P=3}}\label{sec:certificate}
The preceding section determines the degrees and reciprocal symmetries
of $L_{81}$ and $L_{64}$, but not their coefficients.  We now explain
how a finite collection of exact atom counts determines both
polynomials uniquely.  The role of the computation is not to guess a
recurrence: the cohomological dimension and duality arguments already
supply the finite-dimensional constraints that make the reconstruction
rigorous.

\subsection{Exact enumeration of atoms}
For each $H$ and each $a\in\F_{3^H}$, consider the companion product
\[
 M_{a,H}=C_aC_{a^3}\cdots C_{a^{3^{H-1}}}
\]
from Section~\ref{sec:setup}.  By the companion-module criterion,
$f_a$ is an atom if and only if the characteristic polynomial of
$M_{a,H}$ is irreducible over $\F_3$.  Since $M_{a,H}$ is symplectic,
its characteristic polynomial has the self-reciprocal form
\[
 T^4+c_1T^3+c_2T^2+c_1T+1.
\]
The only irreducible self-reciprocal quartics over $\F_3$ are
\[
 T^4+T^3+T^2+T+1,\qquad
 T^4+2T^3+T^2+2T+1.
\]
Thus the trace and second characteristic coefficient suffice to
recognize the two permissible characteristic polynomials.  In
particular, the test can be performed using finite-field arithmetic
alone, without factoring a skew polynomial for every parameter.

For reproducibility, use
$\F_{3^H}=\F_3[X]/(X^H+\sum_{i=0}^{H-1}b_iX^i)$.
The second row below records the integers
$\sum_{i=0}^{H-1} b_i3^i$, with $b_i\in\{0,1,2\}$.
These moduli are irreducible over $\F_3$.
\[
\begin{array}{c|rrrrrrrrrrrr}
H&1&2&3&4&5&6&7&8&9&10&11&12\\\hline
\sum_i b_i3^i&1&1&7&5&7&5&11&11&64&19&11&11\\
A_H(3)&1&3&7&15&50&147&456&1231&3976&11770&35454&106239
\end{array}
\]
Loop over the field elements $a$, form $M_{a,H}$ using the
successive cubes of $a$, and count an element when
\[
 \Tr(M_{a,H})\in\{1,2\},\qquad
 \frac{\Tr(M_{a,H})^2-\Tr(M_{a,H}^2)}2=1 .
\]
Equivalently, one may take one representative per cubing orbit
and weight each accepted representative by its orbit size.
The companion products for consecutive elements of a cubing orbit
are conjugate by cyclic rotation of the invertible factors.

\subsection{From point counts to a determinant quotient}
Put
\[
 b_H=3^H+1-5A_H(3).
\]
By the trace identity \eqref{eq:trace-three},
\begin{align*}
 b_H=\Tr\bigl(F^H\mid H_c^1(\mathbb A^1,\mathcal E_{81})\bigr)-\Tr\bigl(F^H\mid H_c^1(\mathbb A^1,\mathcal E_{64})\bigr).
\end{align*}
Recall the formal trace--determinant identity
\begin{equation}\label{eq:trace-det-reconstruction}
 \det(1-uT)=\exp\left(-\sum_{H\ge1}
 \frac{\Tr(T^H)}{H}u^H\right)
\end{equation}
for any endomorphism $T$ of a finite-dimensional vector space in
characteristic zero.  Consequently the formal series
\begin{equation}\label{eq:reconstruction-series}
 R(u):=\exp\left(-\sum_{H\ge1}\frac{b_H}{H}u^H\right)
 =\frac{L_{81}(u)}{L_{64}(u)}
\end{equation}
is determined by the atom counts.  We write $R(u)=\sum_{h\ge0}r_hu^h$
with $r_0=1$.  Logarithmic differentiation gives
\[
 uR'(u)=-R(u)\sum_{j\ge1}b_ju^j,
\]
and hence the recursion
\begin{equation}\label{eq:reconstruction-recursion}
 hr_h=-\sum_{j=1}^{h}b_jr_{h-j}\qquad(h\ge1).
\end{equation}
For instance, $b_1=-1$ and $b_2=-5$, so that $r_1=1$ and
$r_2=3$.  All calculations take place over $\mathbb Q$; no
floating-point approximation enters the reconstruction.

\subsection{Uniqueness from reciprocal duality}
Write
\[
 L_{81}(u)=\sum_{i=0}^{12}n_iu^i,\qquad
 L_{64}(u)=\sum_{i=0}^{10}d_iu^i,
 \qquad n_0=d_0=1.
\]
Use the relations \eqref{eq:reciprocity-three}.
For $\varepsilon=1$, the unknowns are $n_1,\ldots,n_6$ and
$d_1,\ldots,d_5$. For $\varepsilon=-1$, the middle coefficient
$d_5$ vanishes, leaving ten unknowns.
By \eqref{eq:reconstruction-series},
\[
 L_{64}(u)R(u)=L_{81}(u).
\]
Comparison of coefficients of $u^h$, for $1\le h\le11$, yields
\begin{equation}\label{eq:reconstruction-linear-system}
 \sum_{i=0}^{\min(h,10)}d_i r_{h-i}=n_h,
 \qquad 1\le h\le11,
\end{equation}
where coefficients with indices greater than half the respective
degrees are replaced using reciprocity.  The $r_h$ occurring here
are obtained solely from $A_1(3),\ldots,A_{11}(3)$ by
\eqref{eq:reconstruction-recursion}.
For $\varepsilon=-1$, order the ten unknowns as
$(n_1,\ldots,n_6,d_1,\ldots,d_4)$ and write the equations with
left-hand side $\sum_i d_i r_{h-i}-n_h$.
The $11\times 11$ matrix formed by appending the right-hand
side to the coefficient matrix has determinant $461\,761\,074$,
so that system is inconsistent. Hence $\varepsilon=1$.

For this sign, in the ordering
$(n_1,\ldots,n_6,d_1,\ldots,d_5)$, the coefficient matrix
has determinant $-6\,975\,000$, which is nonzero.  Therefore these
counts determine the two polynomials uniquely.  Solving gives
\[
 (n_1,\ldots,n_6)=(1,2,1,2,17,26),\qquad
 (d_1,\ldots,d_5)=(0,-1,-3,3,18).
\]
The remaining coefficients are fixed by reciprocity, yielding
$L_{81}=N$ and $L_{64}=D$ as in Theorem~\ref{thm:main}.
Substitution into \eqref{eq:reconstruction-series} then proves the
formula for $A_H(3)$ for every $H\ge1$, not merely for the extension
degrees in the table.

The value at $H=12$ is an additional consistency check: it is not
needed for uniqueness, since the eleven coefficients were already
determined from $H\le11$. Substituting the reconstructed
polynomials in \eqref{eq:reconstruction-series} recovers the
twelfth count.

\subsection*{Code and data availability}
The table and the exact reconstruction can be reproduced with the
companion file \texttt{verify\_counts.py}. With NumPy and SymPy
installed, run \texttt{python3 verify\_counts.py}. The defining
field moduli and the count data are also recorded above.

\section{The exact minimal recurrence order for \texorpdfstring{$P=3$}{P=3}}\label{sec:recurrence-three}
\begin{proposition}\label{prop:no-cancel-three}
The polynomials $N$ and $D$ in Theorem~\ref{thm:main}
are coprime and separately squarefree over $\mathbb Q$.
The reduced ordinary generating series for $A_H(3)$ has denominator
\begin{equation}\label{eq:recurrence-denominator-three}
 Q_3(u)=(1-u)(1-3u)N(u)D(u).
\end{equation}
In particular, $A_H(3)$ has a homogeneous constant-coefficient
linear recurrence of \emph{minimal} order $24$ over $\mathbb Q$,
including among recurrences valid only for sufficiently large $H$.
\end{proposition}
\begin{proof}
The Euclidean algorithm over $\F_7$ gives
\begin{equation}\label{eq:gcd-mod-seven}
 \gcd(\overline N,\overline D)=1,
\end{equation}
and over $\F_5$ it gives
\begin{equation}\label{eq:squarefree-mod-five}
 \gcd(\overline N,\overline{N'})=
 \gcd(\overline D,\overline{D'})=1.
\end{equation}
In addition, direct multiplication gives
\begin{align*}
N(u)&=(1+u+3u^2)N_{10}(u),\\
N_{10}(u)&=1-u^2+2u^3+3u^4+8u^5+9u^6+18u^7-27u^8+243u^{10}.
\end{align*}
The polynomial $N_{10}$ is irreducible modulo $17$: after reducing
modulo $N_{10}$ over $\F_{17}$ one has $u^{17^{10}}=u$, and the
gcds of $N_{10}$ with $u^{17^5}-u$ and $u^{17^2}-u$ are both $1$.
Its degree is preserved modulo $17$, so it is irreducible over
$\mathbb Q$.

The leading coefficients of $N$ and $D$, respectively $729$ and $243$, are
nonzero modulo both $5$ and $7$. Thus these modular computations
prove coprimality and squarefreeness over $\mathbb Q$.

The $12$ inverse roots of $N$ and the $10$ inverse roots of $D$
are now distinct within and across the two sets. They all have
absolute value $\sqrt3$; in particular they cannot equal $1$
or $3$. By Theorem~\ref{thm:main}, $5A_H(3)$ is a sum of
$24$ exponential sequences with distinct nonzero bases and
nonzero coefficients. Each contributes a simple pole at the
reciprocal of its base, giving \eqref{eq:recurrence-denominator-three}.
If a polynomial in the shift annihilates $A_H(3)$ for all
sufficiently large $H$, the Vandermonde determinant applied
to $24$ consecutive values forces the polynomial to vanish
at every base. Conversely, the polynomial with precisely these
$24$ roots annihilates the sequence and has rational coefficients
by \eqref{eq:recurrence-denominator-three}.
\end{proof}

\begin{remark}[Explicit recurrence coefficients]
Writing
\[
 Q_3(u)=1+q_1u+\cdots+q_{24}u^{24},
\]
the recurrence takes the form
\[
 A_{H+24}(3)+q_1A_{H+23}(3)+\cdots+q_{24}A_H(3)=0
 \qquad(H\ge1).
\]
The coefficients $q_i$ are obtained by multiplying the four
explicit polynomials in \eqref{eq:recurrence-denominator-three};
no additional counting or cohomological calculation is required.
This expresses all further counts from $24$ consecutive initial
values. The order $24$ is a structural statement, rather than
merely an upper bound obtained from the cohomological degrees.
\end{remark}

\part{The flag curve and the general counting series}
Throughout Part II, assume that $P=p^f$ is odd and $p\nmid f$,
except where a direct finite-field identity is stated for every
odd prime power.
The two factors computed in Part I come from the Steinberg
representation and a regular cuspidal representation. We now
determine their general counterparts and realize the Steinberg
factor as the zeta numerator of a flag curve.

\section{Cohomological degrees and factorization for general \texorpdfstring{$P$}{P}}\label{sec:general}
The conductor data of Lemmas~\ref{lem:inertia} and~\ref{lem:swan}
can be computed from four unipotent character values. Write $r_+,r_-$ for the two rank-one classes
in the Siegel radical and $s_+,s_-$ for its split and anisotropic
rank-two classes. For a representation of dimension $v$ and character
$\chi$, put
\begin{align*}
a_\chi=\frac{1}{P^3}\bigl[&v+\tfrac{P^2-1}{2}
 (\chi(r_+)+\chi(r_-))
 +\tfrac{P(P^2-1)}2\chi(s_+)
 +\tfrac{P(P-1)^2}2\chi(s_-)\bigr],\\
b_\chi=\frac{1}{P^2}\bigl[&v+\tfrac{P^2-1}{2}
 (\chi(s_+)+\chi(s_-))\bigr].
\end{align*}
Then $a_\chi=\dim E^U$, $b_\chi=\dim E^{U_0}$.
To justify evaluating characters on congruence classes of
symmetric forms, note that the Levi subgroup acts by actual
conjugation in $G$:
\[
 \begin{pmatrix}A&0\\0&(A^{\mathsf t})^{-1}\end{pmatrix}
 u_B
 \begin{pmatrix}A^{-1}&0\\0&A^{\mathsf t}\end{pmatrix}
 =u_{ABA^{\mathsf t}}.
\]
Thus any character is constant on each congruence class.
Here $s_+$ and $s_-$ denote the split and anisotropic rank-two
forms, respectively; the split condition is that $-\det B$ be a
square, rather than necessarily that $\det B$ be a square.
To see the second count, identify the two-dimensional summand of symmetric
forms with $x\mapsto\Tr_{\F_{P^2}/\F_P}(b x^2)$, $b\in\F_{P^2}$.
Every nonzero such form is nonsingular, since its polar form is
$(x,y)\mapsto2\Tr_{\F_{P^2}/\F_P}(bxy)$ and the trace pairing
is nondegenerate. Up to a fixed nonzero square class, the
determinant of this bilinear form is
$N_{\F_{P^2}/\F_P}(b)$. Since the norm map is surjective,
its two square classes each occur $(P^2-1)/2$ times.
The first count is the elementary congruence classification
of symmetric two-by-two matrices. In particular, rank-one
matrices have two congruence classes according to the nonzero
square class of their diagonal coefficient; their cardinalities
are added when the character values coincide. Nondegenerate
forms have two classes distinguished by their determinant square
class (or equivalently by split versus anisotropic type,
with the convention just specified).

For nontrivial irreducible characters of $S=\operatorname{PSp}_4(\F_P)$,
the exact degree of the associated compactly supported factor is
\begin{equation}\label{eq:general-degree}
\delta_\chi=
\frac{v-a_\chi}{P-1}+\frac{P}{P+1}(v-b_\chi)-v.
\end{equation}
In particular, this is an exact degree formula, not merely an upper bound.
Recall the cyclic torus $T_0\subset S$ of order
$m=(P^2+1)/2$ introduced in Section~\ref{sec:setup}.
For each irreducible character $\chi$ of $S$, put
$K_\chi=m\dim E_\chi^{T_0}-\dim E_\chi$.
The elementary induction identity
\[
4m\,1_{\rm irr}=m\operatorname{Ind}_{T_0}^S1-\operatorname{Reg}_S
=\sum_\chi K_\chi\chi
\]
then gives the uniform finite factorization
\begin{equation}\label{eq:general-product}
\mathcal R_P(u)=(1-Pu)^{-(m-1)}
 \prod_{\chi\ne1} L_\chi(u)^{K_\chi},\qquad
 \deg L_\chi=\delta_\chi.
\end{equation}
The identity follows because every nonidentity element of $T_0$ is
regular and irreducible, with centralizer $T_0$ and Weyl-group orbit
of size four. Indeed, if a lift $t\in T\setminus\{\pm I\}$ were
conjugate to $-t$, then $t^{P^j}=-t$ for some $j$ by the
description of its eigenvalues in the norm-one torus.
Squaring would give $(t^2)^{P^j}=t^2$. The order of $t^2$
divides $m$, and $\gcd(m,P^2-1)=1$, so for $t^2\ne1$ its
$P$-Frobenius orbit has length four. This forces $4\mid j$,
which would give $t=-t$. Thus the preimage of $C_S(t_0)$ is
$C_G(t)=T$. The same orbit argument gives the four conjugates
of $t_0$ in $T_0$. At the identity the two permutation terms cancel.
Moreover each factor is determined by a finite, explicitly specified
calculation: if
\[S_{\chi,h}=\sum_{a\in\F_{P^h}}\chi(F_{a,h}),\]
then
\[L_\chi(u)=
 \left[\exp\left(\sum_{h=1}^{\delta_\chi}
 S_{\chi,h}\frac{u^h}{h}\right)\right]_{\leq\delta_\chi}.
\]
Thus there are no unknown extension-dependent degrees. We now evaluate
the unipotent contributions and make the remaining family explicit.
The monodromy theorem used here is invoked only under $p\nmid f$;
no assertion based on that theorem is made for $p\mid f$.

\section{The irreducibility indicator and its character constituents}\label{sec:characters}
The Coxeter-torus identity \eqref{eq:indicator-general}, already
proved in Section~\ref{sec:setup}, isolates the unipotent characters
$\eta_\pm$, the Steinberg character, and the regular cuspidal
$\rho_\theta$. We now evaluate their local conductor contributions.
The two additional unipotent characters $\gamma$ and $\xi$ will be
used for permutation representations of lines and flags.

Here are all unipotent values needed to evaluate the conductors. The two
rank-one classes have the same value in each row.
\begin{equation}\label{eq:unipotent-values}
\begin{array}{c|r|rrr}
\chi&\chi(1)&r_\pm&s_+&s_-\\\hline
\eta_-&P(P-1)^2/2&-P(P-1)/2&0&P\\
\eta_+&P(P+1)^2/2&P(P+1)/2&P&0\\
\gamma&P(P^2+1)/2&-P(P-1)/2&P&0\\
\xi&P(P^2+1)/2&P(P+1)/2&0&P\\
\mathrm{St}&P^4&0&0&0\\
\rho_\theta&(P^2-1)^2&1-P^2&1-P&1+P
\end{array}
\end{equation}
The five unipotent rows agree with Srinivasan's original
character table for $\Sp_4(\F_P)$ \cite[\S8]{Srinivasan}, with the
class labels matched by rank and split type. For a complementary
classification in the conformal symplectic group, see Shinoda's
original character table for $\operatorname{CSp}(4,P)$
\cite{Shinoda}; compare also \cite[Table 20]{Breeding} and the
generic unipotent table \cite{Generic}. The latter is based on
Shinoda's work, but we retain Srinivasan's labels for characters
of $\Sp_4(\F_P)$. For the last row, restrictions of
$R_T^G(\theta)$ to unipotent elements are independent of $\theta$.
The Coxeter Green function, equivalently obtained from
\cite[Examples 5.2(a), 5.8]{Geck}, has the values shown.
These also follow from $R_T^G(1)=1-\eta_++\eta_-+\mathrm{St}$.
Substitution in \eqref{eq:general-degree} gives exact, uniform identities:
\begin{equation}\label{eq:general-dimensions}
\begin{array}{c|rrr}
\chi&\dim E_\chi^U&\dim E_\chi^{U_0}&\delta_\chi\\\hline
\eta_-&0&P-1&0\\
\eta_+&P+1&P+1&1\\
\gamma&1&P&1\\
\xi&P&P&0\\
\mathrm{St}&P&P^2&P^2+P\\
\rho_\theta&0&P^2-1&P^2+P-2
\end{array}
\end{equation}
In particular $L_{\eta_-}=L_\xi=1$.

\begin{lemma}\label{lem:linear-factor}
For every $H\geq1$ one has
\[\sum_{a\in\F_{P^H}}\eta_+(F_{a,H})=-1,\qquad L_{\eta_+}(u)=1-u.
\]
\end{lemma}
\begin{proof}
The permutation character on $\F_P$-lines in the natural space is
$1+\eta_++\xi$; this is the Klingen permutation character, see
\cite[Table 15]{Breeding}. A line spanned by $x\ne0$ is parametrized by
$z=x^{P-1}$. Dividing $L_a(x)=0$ by $x$ gives
\[
a=-\frac{z^{n_4}+z^{n_3}+z+1}{z^{P+1}},\qquad
n_j=\frac{P^j-1}{P-1}.
\]
The formula describes the root-line cover by the parameter
$z\in\mathbb G_m$. Indeed, two nonzero roots have the same $z$
precisely when they differ by an element of $\F_P^\times$.
Conversely, if $z\in\F_{P^H}^\times$ and $x^{P-1}=z$, then
\[
 \frac{x^{P^H}}{x}=
 z^{1+P+\cdots+P^{H-1}}
 =N_{\F_{P^H}/\F_P}(z)\in\F_P^\times,
\]
so its root line is Frobenius-stable. Moreover direct differentiation in
characteristic $p$ gives
\[
 \frac{da}{dz}=z^{-(P+2)},
\]
which is nonzero on $\mathbb G_m$, confirming that this description
is an \'etale root-line cover. Consequently its $\F_{P^H}$-points
count all Frobenius-stable root lines over parameters
$a\in\F_{P^H}$, and the total is $P^H-1$.
Since $L_\xi=1$, its character sum is zero for every $H$. The stated
identity follows, and the rank-one cohomological factor is $1-u$.
\end{proof}

The line-cover computation determines one of the degree-one factors.
The remaining geometric input is the compactification of the covers
parametrizing Lagrangian planes and isotropic flags, treated in the
next section; the associated boundary points are counted by inertia
orbits rather than by the character table.

\section{The flag curve and the Steinberg factor}\label{sec:flag}

\subsection{Why full flags appear}
A \emph{full isotropic flag} in a four-dimensional symplectic
$\F_P$-space $V$ is a pair $\ell\subset M$ where
$\ell$ is a one-dimensional subspace and $M$ is a
Lagrangian two-plane. The finite $G$-set of such pairs will be
denoted by $\mathrm{Flags}$, and we likewise write
$\mathrm{Lines}$ and $\mathrm{Lag}$ for the sets of
lines and Lagrangian planes. The latter are finite sets of
$\F_P$-rational subspaces, not varieties whose geometric
points are being counted here.

The Steinberg factor in the preceding character identity has a natural
geometric origin.  The spherical Tits building of
$G=\Sp_4(\F_P)$ is the simplicial complex whose simplices
encode proper nonzero isotropic flags. In rank two it is a
bipartite graph: its vertices are the isotropic
lines and the Lagrangian planes in a four-dimensional symplectic space,
and its edges are the incident pairs $\ell\subset M$.  Its first
reduced homology is the Steinberg representation.  In particular,
the augmented chain complex gives an exact sequence of $G$-modules
\begin{equation}\label{eq:flag-building-complex}
 0\longrightarrow \mathrm{St}\longrightarrow
 \overline{\mathbb Q}_\ell[\mathrm{Flags}]
 \longrightarrow
 \overline{\mathbb Q}_\ell[\mathrm{Lines}]
 \oplus\overline{\mathbb Q}_\ell[\mathrm{Lag}]
 \longrightarrow\overline{\mathbb Q}_\ell\longrightarrow0.
\end{equation}
Here the last map takes the difference of the two augmentation maps.
Consequently the corresponding permutation characters satisfy
\begin{equation}\label{eq:flag-building-character}
 \mathbf1_{\mathrm{Flags}}=
 \mathbf1_{\mathrm{Lines}}+\mathbf1_{\mathrm{Lag}}-1+\mathrm{St}.
\end{equation}
For the root spaces $V_a$ of Section~\ref{sec:setup}, this suggests
studying the finite cover of the parameter line that classifies
\emph{isotropic full flags} $\ell\subset M\subset V_a$.
Although the flags form part of a positive-dimensional flag variety
before taking $\F_P$-points, their set in each $V_a$ is finite;
therefore the resulting cover of the one-dimensional parameter line
is itself a curve.

We shall write down this flag cover explicitly and show that its
smooth projective model realizes the Steinberg $L$-polynomial.
The lower-level line and Lagrangian-plane covers account for the
other constituents of~\eqref{eq:flag-building-character}.
The same building-theoretic construction in higher symplectic rank is
explained in Appendix~\ref{app:higher-flags}; the special feature
here is that a \emph{single} smooth projective curve suffices.

\subsection{Inertia orbits and the boundary}
The inertia subgroup computed in Section~\ref{sec:monodromy}
controls how these finite covers compactify.  Indeed, geometric
points over infinity in the normalization of a finite cover of the
affine line correspond to inertia orbits in its geometric generic
fibre.  The following calculation therefore determines the number
of boundary points; the assertion about rationality uses the
arithmetic decomposition group, not just geometric inertia.

\begin{lemma}[Inertia orbits on isotropic flags]\label{lem:flag-orbits}
Let $V=W\oplus W^*$ be symplectic of dimension four, and let
$I=U\rtimes \Sigma$, where $U$ is the Siegel radical of the stabilizer of
$W$ and $\Sigma\simeq\F_{P^2}^{\times}$ acts on $W$ as a Singer group.
Then $I$ has three orbits on Lagrangian planes $M$, indexed by
$\dim(M\cap W)=2,1,0$, and four orbits on flags $\ell\subset M$
with $\dim\ell=1$ and $M$ Lagrangian. The plane orbit sizes are
$1$, $P(P+1)$, and $P^3$, respectively. The flag orbits are
\[
M=W;\quad \dim(M\cap W)=1,\ \ell=M\cap W;\quad
\dim(M\cap W)=1,\ \ell\ne M\cap W;\quad M\cap W=0.
\]
Their sizes in this order are $P+1$, $P(P+1)$, $P^2(P+1)$,
and $P^3(P+1)$. For the root-space family, the positive-valuation plane $W$ is
preserved by the arithmetic decomposition group at infinity.
Consequently these orbits are Frobenius-stable and correspond to rational
boundary points in the smooth projective compactifications of the
respective connected covers.
\end{lemma}
\begin{proof}
The Singer subgroup is transitive on the lines in $W$ and $W^*$.
The Lagrangian complements to $W$ are graphs of symmetric maps
$W^*\to W$, on which $U$ acts simply transitively; the stabilizer
of one complement contains $\Sigma$, transitive on its lines. This proves
the orbit assertion when $M\cap W=0$, and $M=W$ is immediate.
For $\dim(M\cap W)=1$, fix $L=M\cap W$. The planes $M$ correspond
to lines in the symplectic plane $L^\perp/L$ other than $W/L$.
The unipotent radical acts transitively on these $P$ lines.
Within such an $M$, the stabilizer of $M$ in $U$ is transitive on
the lines $\ell\ne L$. To see this, write $M=L+\F_P(w+\lambda)$,
where $\lambda\in W^*$ is nonzero and annihilates $L$.
For any $h\in L$, a symmetric map $W^*\to W$ can send
$\lambda$ to $h$. Indeed, in a basis of $W^*$ starting with
$\lambda$, the first column of its symmetric matrix can be
prescribed arbitrarily. The corresponding element of $U$ sends
$w+\lambda$ to $w+\lambda+h$, giving all translations on the
$P$ lines of $M$ other than $L$. Since $\Sigma$ is transitive on possible $L$,
there are precisely the two stated flag orbits in this case.
Finally the boundary geometric points of a finite \'etale cover of
$\mathbb A^1$ are the inertia orbits on its generic fiber.
The intersection conditions with $W$ are preserved by the whole
decomposition group, so each boundary orbit is individually stable
under Frobenius and hence rational over $\F_P$.
\end{proof}

\subsection{An explicit equation and its cohomology}
We now express the full-flag moduli problem by two coordinates
and compare the cohomology of its compactification with the
Steinberg factor.

\begin{theorem}\label{thm:flagcurve}
Assume $P=p^f$ is odd and $p\nmid f$. Let $\mathcal B_P$ be the smooth
projective model of the affine curve
\begin{equation}\label{eq:flag-curve}
z^{P^2+1}s^{P+1}+(1+z+z^{P^2+1})s+1=0,
\qquad z,s\ne0.
\end{equation}
Then $\mathcal B_P$ is geometrically connected, has genus $P(P+1)/2$,
and the complement in $\mathcal B_P$ of the displayed affine curve consists of four
$\F_P$-rational points. Its zeta numerator is exactly $L_{\mathrm{St}}$:
\[
Z(\mathcal B_P,u)=\frac{L_{\mathrm{St}}(u)}{(1-u)(1-Pu)}.
\]
The compactification of the Lagrangian-plane cover is a genus-zero curve
whose affine part has three rational points removed. In particular,
\[\sum_{a\in\F_{P^H}}\#\{\text{Frobenius-stable Lagrangian planes in }V_a\}=P^H-2
\qquad(H\ge1).\]
\end{theorem}
\begin{proof}
\emph{The flag coordinates.}
The coordinates below have an intrinsic meaning: $z$ records the
root line $\ell$, while $s$ records the additional choice of
Lagrangian plane $M$ containing it. More explicitly, for a nonzero
root $x$ and $z=x^{P-1}$, the equation $L_a(x)=0$ implies
\begin{equation}\label{eq:flag-line-parameter}
 a=-\frac{z^{1+P+P^2+P^3}+z^{1+P+P^2}+z+1}
          {z^{1+P}}.
\end{equation}
Thus the $z$-coordinate determines both the line and its
parameter, and the two-coordinate equation below is to be
viewed as a finite cover of the $a$-line through
\eqref{eq:flag-line-parameter}. Both are invariant under
changing the chosen generators of the flag. The resulting plane
is required to lie in the root space $V_a$ and to be isotropic
for the pairing $B$.
Using the form $B$ of Section~\ref{sec:setup}, fix a root line,
choose $x\ne0$ on it, and put
\[
 y=xv,\qquad w=v^P-v,\qquad z=x^{P-1}.
\]
Dividing the expression for $B(x,xv)$ by $-x^{P^2+1}$ gives
\begin{align*}
 -\frac{B(x,xv)}{x^{P^2+1}}
 &= (v^{P^2}-v)+z(v^{P^2}-v^P)
       +z^{P^2+1}(v^{P^3}-v^P)\\
 &=w+w^P+zw^P+z^{P^2+1}(w^{P^2}+w^P).
\end{align*}
Thus $B(x,y)=0$ is equivalent to
\[
z^{P^2+1}w^{P^2}+(1+z+z^{P^2+1})w^P+w=0.
\]
The nonvanishing of $w$ expresses the linear independence of $x,y$.
Putting $s=w^{P-1}$ and dividing the equation by $w$ gives
\eqref{eq:flag-curve}, since $w^{P^2-1}=s^{P+1}$.
Changing the basis of the flag preserves $z,s$: replacing
$x$ by $\lambda x$ with $\lambda\in\F_P^\times$ fixes $z$, and
replacing $y$ by $\mu y+\nu x$ with $\mu\in\F_P^\times$
sends $w=(y/x)^P-y/x$ to $(\mu/\lambda)w$,
leaving $w^{P-1}$ fixed.
Conversely, $z$ specifies the root line, while $s$ specifies the
$\F_P$-line spanned by $w$ in the Artin--Schreier quotient;
choosing a lift $v$ with $v^P-v=w$ reconstructs the two-plane
$\langle x,xv\rangle_{\F_P}$. Two lifts $v$ of the same $w$
differ by an element of $\F_P$, so the corresponding roots $xv$
differ by an element of $\F_P x$ and give the same plane.
The coordinate map and its inverse on geometric points commute with
$F_Q$: applying $Q$-power Frobenius to $x,y$ sends the coordinates
to $(z^Q,s^Q)$. By the bijection, a flag is $F_Q$-stable precisely
when both $z$ and $s$ belong to $\F_Q$.
The alternating-form equation is exactly the displayed equation in
$w$; the identity
\[
B(x,y)^P-B(x,y)=L_a(x)y^{P^2}-x^{P^2}L_a(y)
\]
shows that $B(x,y)=0$ and $L_a(x)=0$ force $L_a(y)=0$.
Thus \eqref{eq:flag-curve} describes the full-flag cover.

\emph{Smoothness and connectedness.}
If $\Phi(z,s)$ denotes the left-hand side of
\eqref{eq:flag-curve}, then on the curve
\[
 \frac{\partial\Phi}{\partial s}
 =z^{P^2+1}s^P+1+z+z^{P^2+1}=-s^{-1}.
\]
Thus the projection to the root-line coordinate $z$ is \'etale
on the displayed open set. Together with the root-line cover
this gives the requisite finite \'etale full-flag cover over
the parameter line. Full geometric monodromy, which acts
transitively on full flags, proves geometric connectedness.

\emph{The Lagrangian-plane cover.}
By Lemma~\ref{lem:flag-orbits}, the plane cover has three
$\F_P$-rational boundary points. Its permutation character is
$1+\eta_++\gamma$, and the dimensions of the corresponding
affine cohomology groups in~\eqref{eq:general-dimensions}
give $\dim H_c^1=2$.  For a connected smooth affine curve
obtained by removing $r$ geometric points from a projective
curve of genus $g$, one has
\[
 \dim H_c^1=2g+r-1.
\]
Hence $2=2g+3-1$, and $g=0$. Since all three boundary points are
rational, the affine plane cover has $P^H-2$ rational points
over $\F_{P^H}$. The resulting character sum gives
$L_\gamma(u)=1-u$, alongside the earlier
$L_{\eta_+}(u)=1-u$ and $L_\xi(u)=1$.

\emph{The full-flag cover.}
The flag permutation character is
\[
 1+2\eta_++\gamma+\xi+\mathrm{St},
\]
in agreement with~\eqref{eq:flag-building-character} and the
line and plane decompositions.  Let $Y$ be its affine cover.
The preceding linear factors and the permutation-sheaf decomposition
yield the Frobenius characteristic polynomial
\begin{equation}\label{eq:flag-affine-det}
 \det(1-uF\mid H_c^1(Y_{\overline{\F}_P}))
 =(1-u)^3L_{\mathrm{St}}(u).
\end{equation}
In particular its first compactly supported cohomology has
dimension $P^2+P+3$. The Steinberg factor has degree
$P^2+P$ by~\eqref{eq:general-dimensions}.

\emph{Removing the boundary contribution.}
Let $D=\mathcal B_P\setminus Y$. By the inertia-orbit lemma,
$D$ consists of four $\F_P$-rational points. The localization
sequence gives a Frobenius-equivariant exact sequence
\[
0\longrightarrow H^0(\mathcal B_P)
\longrightarrow H^0(D)
\longrightarrow H_c^1(Y)
\longrightarrow H^1(\mathcal B_P)\longrightarrow0,
\]
where all cohomology is taken after base change to
$\overline{\F}_P$ with coefficients in $\overline{\mathbb Q}_\ell$.
The cokernel of $H^0(\mathcal B_P)\to H^0(D)$ is
three-dimensional with trivial Frobenius action. Therefore
\[
 \det(1-uF\mid H_c^1(Y))
 =(1-u)^3\det(1-uF\mid H^1(\mathcal B_P)).
\]
Comparing with~\eqref{eq:flag-affine-det} identifies the latter
determinant with $L_{\mathrm{St}}(u)$.
Its degree is $P^2+P=2g(\mathcal B_P)$, proving the genus formula
and the asserted zeta function.
\end{proof}

\begin{remark}[A direct low-degree check at $P=3$]
The affine equation \eqref{eq:flag-curve}, with $z,s\ne0$,
can be counted independently over small fields. For
$H=1,2,3,4,5$, its numbers of points over $\F_{3^H}$ are
\[
 1,\quad 9,\quad 22,\quad 81,\quad 321.
\]
Adding the four rational boundary points gives
$5,13,26,85,325$, respectively. These agree with
$3^H+1-\sum_i\alpha_i^H$ calculated from the explicit
polynomial $N(u)$ of Theorem~\ref{thm:main}. This is a
consistency check, not a substitute for the geometric proof of
the flag-curve identification.
\end{remark}

\section{The counting series for general \texorpdfstring{$P$}{P}}\label{sec:uniform}
The character identity of Proposition~\ref{prop:generic-indicator}
reduces the atom-counting problem to the Steinberg and regular
cuspidal constituents. The point of this section is to assemble
these contributions into one formula, while keeping distinct the
geometrically realized factor and the factors not yet explicitly
evaluated for arbitrary $P$.

\subsection{The two cohomological factors}
Throughout this section let $P=p^f$ be odd, with $p\nmid f$, and put
\[
 m=\frac{P^2+1}{2},\qquad \kappa=\frac{P^2-1}{8}.
\]
Recall the set $\mathcal O$ of Frobenius orbits of nontrivial
characters of the Coxeter torus $T_0$ from
Section~\ref{sec:setup}. Set
\begin{equation}\label{eq:two-factors}
 S_P(u)=L_{\mathrm{St}}(u),\qquad
 C_P(u)=\prod_{\theta\in\mathcal O}L_{\rho_\theta}(u).
\end{equation}
Here $L_\chi(u)$ denotes the Frobenius characteristic polynomial
on $H_c^1(\mathbb A^1_{\overline{\F}_P},\mathcal E_\chi)$, with
constant term one. The conductor computation in
Sections~\ref{sec:general}--\ref{sec:characters} gives
\begin{equation}\label{eq:two-factor-degrees}
 \deg S_P=P^2+P,\qquad
 \deg C_P=\kappa(P^2+P-2).
\end{equation}
Both polynomials have integral coefficients. For $S_P$ this follows
also from Theorem~\ref{thm:flagcurve}, since it is the zeta numerator
of the curve $\mathcal B_P$. For $C_P$, the roots are algebraic
integers and the full collection of regular cuspidal constituents is
stable under Galois conjugation; hence the product has rational
algebraic-integer coefficients. In the present normalization their
inverse roots have complex absolute value $\sqrt P$. In particular,
none of these inverse roots equals $1$ or $P$.

\begin{remark}[Purity of the two factors]\label{rem:purity-factors}
Here purity means that every inverse root $\alpha$ of the relevant
Frobenius polynomial satisfies
$|\iota(\alpha)|=P^{1/2}$ for every embedding
$\iota:\overline{\mathbb Q}\hookrightarrow\mathbb C$.
For $S_P$, this follows from
Theorem~\ref{thm:flagcurve}: it is the characteristic polynomial
of Frobenius on $H^1$ of the smooth projective curve $\mathcal B_P$,
which is pure of weight one by the Weil conjectures.
For each regular cuspidal factor $L_{\rho_\theta}$, the
corresponding local system has finite monodromy and is therefore
pointwise pure of weight zero. Moreover,
$E_{\rho_\theta}^{I_\infty}=0$, since
$E_{\rho_\theta}^{U}=0$ by
\eqref{eq:general-dimensions} and $U\subset I_\infty$.
To see why this yields purity of its compactly supported
cohomology, trivialize the local system on a finite Galois cover
of $\mathbb A^1$ and take a smooth projective compactification
$X$ of that cover. The boundary contribution to $H_c^1$
is built from inertia-invariant stalks; it vanishes in the
$\rho_\theta$-isotypic component by the preceding equality.
Consequently $H_c^1(\mathbb A^1_{\overline{\F}_P},
\mathcal E_{\rho_\theta})$ is a Frobenius-stable direct
summand of $H^1(X_{\overline{\F}_P},\overline{\mathbb Q}_\ell)$.
This proves purity of weight one for the regular cuspidal factors
as well. The vanishing of inertia invariants is essential here:
finite monodromy alone would in general give only mixed weights
for cohomology of an open curve.
\end{remark}

For a polynomial $E(u)=\prod_i(1-\lambda_i u)$, write
$p_H(E)=\sum_i\lambda_i^H$. We use the notation
\begin{equation}\label{eq:s-c-powers}
 s_H=p_H(S_P),\qquad c_H=p_H(C_P).
\end{equation}
These are algebraic integers and, as the coefficients of the two
polynomials are integers, in fact ordinary integers.

\subsection{The general counting formula}
\begin{theorem}\label{thm:uniform-count}
For every $H\ge1$, the number $A_H(P)$ of atoms in the
one-parameter family satisfies
\begin{equation}\label{eq:uniform-count}
 m A_H(P)=\kappa(P^H+1-s_H)+c_H.
\end{equation}
Equivalently, in terms of the flag curve of Theorem~\ref{thm:flagcurve},
\begin{equation}\label{eq:uniform-geometric-count}
 m A_H(P)=\kappa\#\mathcal B_P(\F_{P^H})+c_H.
\end{equation}
The exponential counting series defined in Section~\ref{sec:setup}
therefore has the rational expression
\begin{equation}\label{eq:uniform-product}
 \mathcal R_P(u)=
 \left[\frac{S_P(u)^\kappa}
 {(1-u)^\kappa(1-Pu)^\kappa C_P(u)}\right]^4.
\end{equation}
These identities are understood as identities of rational functions,
or equivalently of formal series at $u=0$.
\end{theorem}
\begin{proof}
Evaluate the character identity \eqref{eq:indicator-general}
on the Frobenius elements attached to parameters
$a\in\F_{P^H}$ and sum over $a$.
The trivial character contributes $P^H$.
The trace formulas established earlier give, in the same
normalization,
\[
 \sum_a\eta_+(F_{a,H})=-1,\qquad
 \sum_a\eta_-(F_{a,H})=0,\qquad
 \sum_a\mathrm{St}(F_{a,H})=-s_H,
\]
while the sum of the $\rho_\theta$-contributions is $-c_H$.
The signs in the indicator identity thus yield
\[
 m A_H(P)=\kappa(P^H-(-1)+0-s_H)-(-c_H),
\]
which is \eqref{eq:uniform-count}.
By Theorem~\ref{thm:flagcurve},
$\#\mathcal B_P(\F_{P^H})=P^H+1-s_H$, giving
\eqref{eq:uniform-geometric-count}.
Finally, for any $E(u)=\prod_i(1-\lambda_i u)$,
\[
 -\log E(u)=\sum_{H\ge1}\frac{p_H(E)}H u^H.
\]
Multiply \eqref{eq:uniform-count} by $4u^H/H$, sum over $H\ge1$,
and exponentiate, using the definition of $\mathcal R_P$.
This proves \eqref{eq:uniform-product}.
\end{proof}

\begin{remark}[What is and is not explicit]\label{rem:explicit-uniform}
Equation~\eqref{eq:uniform-geometric-count} splits the atom count
into a term given by rational points on an explicit curve and a
regular cuspidal trace term. It is not, by itself, a closed numerical
formula for all $P$, since the coefficients of $C_P$ have not been
fully evaluated. For $P=3$, $\kappa=1$, $S_3=N$, and $C_3=D$;
thus Theorem~\ref{thm:uniform-count} specializes to the completely
explicit formula of Theorem~\ref{thm:main}. For general $P$,
Section~\ref{sec:coefficients} evaluates the initial coefficients of
the two factors. The degrees in \eqref{eq:two-factor-degrees} are the
degrees \emph{before} possible common factors of $S_P$ and $C_P$
are cancelled.
\end{remark}

\section{Poles and minimal recurrences}\label{sec:spectrum}
The rational expression in Theorem~\ref{thm:uniform-count} does more
than encode the counts: it implies that the sequence $A_H(P)$
satisfies a homogeneous linear recurrence with constant coefficients.
A precise minimal order, however, depends on which Frobenius
inverse roots actually survive in the \emph{virtual} quotient of
the Steinberg and cuspidal factors.

\subsection{The ordinary generating series and its poles}
Set
\[
 \mathcal Q_P(u)=\frac{S_P(u)^\kappa}{C_P(u)}.
\]
Taking the logarithmic derivative of \eqref{eq:uniform-product}
gives the \emph{ordinary}, rather than exponential, counting series:
\begin{equation}\label{eq:ordinary-series}
 m\sum_{H\ge1}A_H(P)u^H
 =\frac{\kappa Pu}{1-Pu}+\frac{\kappa u}{1-u}
 +u\frac{\mathcal Q_P'(u)}{\mathcal Q_P(u)}.
\end{equation}
Indeed,
$u\mathcal R_P'(u)/\mathcal R_P(u)=4m\sum_{H\ge1}A_H(P)u^H$.
After cancellation, write
\begin{equation}\label{eq:surviving-spectrum}
 \mathcal Q_P(u)=\prod_{\alpha\in\Sigma_P}(1-\alpha u)^{e_\alpha},
 \qquad e_\alpha\in\mathbb Z\setminus\{0\},
\end{equation}
where the numbers $\alpha$ are distinct. Then
\begin{equation}\label{eq:power-sum-recurrence}
 mA_H(P)=\kappa P^H+\kappa-\sum_{\alpha\in\Sigma_P}
 e_\alpha\alpha^H.
\end{equation}
Thus each surviving inverse root gives a \emph{simple} pole of the
ordinary generating series at $u=\alpha^{-1}$, regardless of its
multiplicity in $\mathcal Q_P$. Multiplicity changes the residue, not the
pole order. A common inverse root cancels only when its signed
multiplicity in \eqref{eq:surviving-spectrum} is zero.

\begin{proposition}[Minimal-order principle]\label{prop:recurrence-principle}
Let
$r_P=|\Sigma_P|$ be the number of \emph{distinct} inverse roots
surviving in the reduced quotient $\mathcal Q_P$.
Then $A_H(P)$ satisfies a homogeneous rational
constant-coefficient recurrence of minimal order
\[
 2+r_P.
\]
Its characteristic polynomial can be taken to be
\begin{equation}\label{eq:general-recurrence-char}
 (T-1)(T-P)\prod_{\alpha\in\Sigma_P}(T-\alpha).
\end{equation}
This order is minimal even if a recurrence is required only for
all sufficiently large $H$.
\end{proposition}
\begin{proof}
By the purity argument in Remark~\ref{rem:purity-factors},
all $\alpha\in\Sigma_P$ satisfy
$|\iota(\alpha)|=\sqrt P$ for every complex embedding
$\iota$; in particular, none equals $1$ or $P$.
The coefficients of $P^H$, $1^H$, and each $\alpha^H$ in
\eqref{eq:power-sum-recurrence} are nonzero, and all bases are
distinct. Their powers are linearly independent over characteristic
zero by the Vandermonde determinant. Hence any eventual
constant-coefficient recurrence must have each of these bases as
a root of its characteristic polynomial. Conversely, the polynomial
\eqref{eq:general-recurrence-char} annihilates every summand.
It has rational coefficients: $\mathcal Q_P\in\mathbb Q(u)$, so the
set of its poles and zeros, including the nonzero signed
multiplicities, is stable under Galois conjugation.
\end{proof}

The indicator identity has already eliminated the elementary
non-Steinberg factors. For general $P$, it does not rule out
cancellation between $S_P$ and $C_P$ or repeated inverse roots
within either factor, and we make no uniform assertion of
minimality without checking those possibilities.

\section{The first two coefficients for general \texorpdfstring{$P$}{P}}\label{sec:coefficients}
We now turn from the general cohomological description of the counting
series to explicit coefficient calculations. There are two rather
different phenomena. Over $\F_P$, the equation of the flag curve
specializes to a quadratic equation whose discriminant defines an
auxiliary genus-one curve; this determines the \emph{linear}
coefficients by the trace of a fixed elliptic curve. Over
$\F_{P^2}$, trace--norm coordinates instead reduce the point counts
to explicit two-variable quadratic-character sums. An identity
between these sums permits a common formula for both \emph{quadratic}
coefficients.

The cohomological factors $S_P$ and $C_P$ are used under the standing
hypothesis $P=p^f$ with $p\nmid f$. The direct arithmetic counts
$A_1(P),A_2(P)$ and the finite-field character-sum identities below
make sense for every odd prime power $P$, without that monodromy
hypothesis. Throughout, $\chi$ denotes the quadratic character of
$\F_P$, extended by $\chi(0)=0$.

\subsection{The linear coefficients: an auxiliary elliptic curve}
Write
\[
S_P(u)=1+a_1u+a_2u^2+\cdots,\qquad
C_P(u)=1+b_1u+b_2u^2+\cdots.
\]
The first coefficient of a Frobenius polynomial is minus its first
power trace. The point-count calculation below makes this trace
accessible without computing a character table.

\begin{theorem}\label{thm:first-coeff}
For $p\ne3,5$, let
\[
\tau_P=P+1-\#E(\F_P),\qquad E:\ y^2=x(x-16)(x+9).
\]
Then
\begin{equation}\label{eq:first-good}
a_1=-\tau_P,\qquad b_1=-\frac{P+1-\tau_P}{4}.
\end{equation}
For $P=3^f$ put $\tau_P=(-1)^f$, and for $P=5^f$ put $\tau_P=1$.
In these two cases,
\begin{equation}\label{eq:first-exceptional}
a_1=-\tau_P,\qquad b_1=\kappa-\frac{P-\tau_P}{4}.
\end{equation}
For good primes, the sequence $\tau_P$ can also be evaluated from the
weight-two newform
\[
\eta(z)\eta(3z)\eta(5z)\eta(15z)=\sum_{n\geq1}e_n q^n.
\]
Namely $\tau_p=e_p$, and for $f\geq2$,
$\tau_{p^f}=e_p\tau_{p^{f-1}}-p\tau_{p^{f-2}}$, with $\tau_{p^0}=2$.
These power traces are not, in general, the newform coefficients $e_{p^f}$.
\end{theorem}
\begin{proof}
At $\F_P$-points, \eqref{eq:flag-curve} reduces to a quadratic in $s$:
\[
z^2s^2+(z^2+z+1)s+1=0.
\]
Its discriminant is
$\Delta(z)=(z^2-z+1)(z^2+3z+1)$.
Adding the four boundary points gives
\begin{equation}\label{eq:flag-first}
\#\mathcal B_P(\F_P)=P+2+\sum_{z\in\F_P}\chi(\Delta(z)),\qquad
a_1=1+\sum_z\chi(\Delta(z)).
\end{equation}
If $p\ne3,5$, the quartic $\Delta$ is squarefree, and the
smooth projective model of $F:w^2=\Delta(z)$ has genus one and two
rational points at infinity. Consequently its point count equals
that of $\mathcal B_P$ over $\F_P$; this is a point-count
identity, not a morphism from the full flag curve to $F$.
We next show that $F$ is isogenous to $E$. The involution
$(z,w)\mapsto(1/z,-w/z^2)$ is fixed-point-free when $p\ne3,5$:
a fixed point would require $z=\pm1$ and $w=0$, which can only
happen in the excluded characteristics. Its quotient has equation
\[
r=z+z^{-1},\quad v=(z-z^{-1})w/z,\qquad
v^2=(r^2-4)(r-1)(r+3).
\]
The substitution
$X=5(r+2)/(r-2)$, $Y=v(X-5)^2/20$
identifies the quotient with $E':Y^2=X(X-1)(X+15)$.
On the other hand, translating $x$ on $E$ by $u=x+9$ gives
$y^2=u^3-34u^2+225u$.
The standard two-isogeny with kernel generated by $(0,0)$ has target
$\eta^2=\xi^3+68\xi^2+256\xi$; writing $\xi=4(X-1)$ and $\eta=8Y$
gives exactly $E'$. Thus $E$ and $F$ have identical Frobenius traces,
proving the formula for $a_1$.

We also obtain the atom count $A_1(P)$ directly, independently
of the flag cover. Because $H=1$, the automorphism $\sigma$ of
$\F_P$ is the identity, so the skew polynomial becomes an ordinary
reciprocal quartic. Substitute $Y=t+t^{-1}$ in this quartic.
The necessary and sufficient conditions for irreducibility are
\begin{equation}\label{eq:reciprocal-test}
\chi(9-4a)=-1,\qquad \chi(a(a+4))=-1.
\end{equation}
Indeed the first condition makes $Y^2+Y+a-2$ irreducible over $\F_P$;
the second says that the norm of $Y^2-4$ is a nonsquare, so the two
remaining roots have degree four. When $p\ne3,5$, expansion of \eqref{eq:reciprocal-test}
has no zero-discriminant correction and gives
\[
A_1(P)=\frac{P+1+\sum_a\chi((9-4a)a(a+4))}{4}
       =\frac{\#E(\F_P)}4.
\]
Under $x=-4a$, the cubic $x(x-16)(x+9)$ becomes
$16a(a+4)(9-4a)$. Since $16=4^2$, the quadratic character of
the cubic equals that of $a(a+4)(9-4a)$.
The coefficient identity from \eqref{eq:uniform-count} is
\[b_1=\kappa\#\mathcal B_P(\F_P)-m A_1(P).\]
Since $m=4\kappa+1$, this proves \eqref{eq:first-good}.

For $p=3$, $\Delta(z)=(z+1)^2(z^2+1)$; for $p=5$,
$\Delta(z)=(z-1)^2(z^2-z+1)$. These factorizations explain
why the genus-one argument above must not be applied in these two
characteristics.
The quadratic character sums in \eqref{eq:flag-first} give $a_1=-(-1)^f$ and $-1$,
respectively. In \eqref{eq:reciprocal-test} exactly one common zero requires correction;
direct evaluation gives $A_1=(P-\tau_P)/4$. This proves \eqref{eq:first-exceptional}.

Finally consider the integral model
\[
 E_{\mathrm{min}}:\quad y^2+xy+y=x^3+x^2-10x-10.
\]
The change of coordinates
\[
 (x,y)\longmapsto (4x+4,\,8y+4x+4)
\]
\emph{from $E_{\mathrm{min}}$ to $E$} gives
$Y^2=X(X-16)(X+9)$, as is checked by direct substitution.
This model belongs to isogeny class $15a$, with the newform stated
above; see \cite{LMFDB15}.
The prime-power recursion is the characteristic-polynomial identity for
its two Frobenius eigenvalues.
\end{proof}

\subsection{The quadratic coefficients: trace--norm descent}
For extension degree $H=2$ the preceding one-variable elliptic
reduction is no longer available in the same form. The useful
replacement is the trace--norm map
\[
  s\longmapsto (T,N)=(s+s^P,s^{P+1})
\]
from $\F_{P^2}$ to $\F_P^2$. The polynomial
$X^2-TX+N$ has discriminant $T^2-4N$; the number of
$s\in\F_{P^2}$ with the given trace and norm is
$1-\chi(T^2-4N)$. This multiplicity will be built into both
character sums.

Define the following sums over the base field. For $T,N\in\F_P$ set
\[
W(T,N)=1-\chi(T^2-4N),\qquad B(T,N)=3T+4N+4,
\]
and put
\begin{align*}
U_P&=\sum_{T,N\in\F_P}W(T,N)
\chi\!\left(N\bigl(B(T,N)^2+TB(T,N)+N\bigr)\right),\\
V_P&=\sum_{T,N\in\F_P}W(T,N)
\nu\bigl((T-1)^2-4N+8\bigr)
\nu\bigl((N+2)^2-4(T-1)^2\bigr),
\end{align*}
where $\nu(x)=(\chi(x)^2-\chi(x))/2$ is exactly the indicator of
nonzero nonsquares. All constants in these sums are interpreted in
$\F_P$; there is no restriction excluding $p=3,5$.

\begin{theorem}\label{thm:second-coeff}
The second coefficients are
\begin{equation}\label{eq:second-coefficients}
a_2=\frac{\tau_P^2+U_P}{2},\qquad
b_2=\frac{b_1^2-m V_P+\kappa(P^2+1+U_P)}2.
\end{equation}
In particular these formulas contain no unevaluated finite-group
characters or unknown Frobenius eigenvalues; the two remaining sums
are explicit quadratic-character sums over the base field.
\end{theorem}
\begin{proof}
Over $\F_{P^2}$, regard \eqref{eq:flag-curve} as a quadratic in $z$:
\[(s^{P+1}+s)z^2+sz+s+1=0.
\]
The exceptional values are $s=0,-1$. The first gives no affine points,
and the second gives only the excluded $z=0$.
For every other $s$, the constant and leading coefficients are nonzero,
so its number of solutions is $1+\chi_2(D(s))$, where
\[
D(s)=s^2-4(s^{P+1}+s)(s+1).
\]
Since $D(0)=0$ and $D(-1)=1$, the projective point count is
$P^2+1+\sum_s\chi_2(D(s))$.
Writing $T=s+s^P$ and $N=s^{P+1}$ gives
\[
D(s)=-(3T+4N+4)s-N,
\quad\operatorname{Norm}(D(s))
=N\bigl((3T+4N+4)^2+T(3T+4N+4)+N\bigr).
\]
There are $W(T,N)$ elements with this trace and norm, and
$\chi_2(x)=\chi(\operatorname{Norm}(x))$.
Thus $\#\mathcal B_P(\F_{P^2})=P^2+1+U_P$, and the second
Steinberg power trace is $-U_P$.

For the atom count at $H=2$, direct multiplication of the two companion
matrices shows that, for $T=a+a^P$ and $N=a^{P+1}$,
\[
\det(XI_4-C_aC_{a^P})
=X^4+(T-1)X^3+NX^2+(T-1)X+1.
\]
The same reciprocal-quartic test used in \eqref{eq:reciprocal-test} gives precisely the two
nonsquare indicators defining $V_P$. The multiplicity of $(T,N)$ is
again $W(T,N)$, so $A_2(P)=V_P$.
Equation \eqref{eq:uniform-count} identifies the second cuspidal power trace as
$mV_P-\kappa(P^2+1+U_P)$.
Newton's identity $[u^2]\prod_i(1-\alpha_i u)
=((\sum_i\alpha_i)^2-\sum_i\alpha_i^2)/2$ proves \eqref{eq:second-coefficients}.
\end{proof}

For example, exact evaluations give
\[
\begin{array}{c|rrrr}
P&a_1&a_2&b_1&b_2\\\hline
3&1&2&0&-1\\5&-1&-2&2&1\\7&0&0&-2&2\\
9&-1&-4&8&28\\11&4&8&-4&8\\13&2&2&-4&8\\
17&-2&9&-4&-30\\19&-4&19&-4&-40\\
25&-1&-12&72&2556\\49&14&49&-16&-160
\end{array}
\]
Theorem~\ref{thm:first-coeff} evaluates the linear coefficients through
a fixed elliptic curve (or its explicit modular form).
Theorem~\ref{thm:second-coeff} gives exact base-field character sums for
the quadratic coefficients. A uniform modular or other closed evaluation
of these sums, and evaluations of all higher coefficients for general
$P$, are not asserted.

\subsection{Simplifying the degree-two atom count}
The sum $V_P$ is an explicit answer, but it involves two nonsquare
indicators. To simplify it, we expand these indicators and evaluate
the six lower-order character sums separately. The only remaining
three-factor sum is $J_P$ below. This calculation is valid over every
finite field of odd order, without the monodromy restriction. Set
\[
K=T^2-4N,\quad D=(T-1)^2-4N+8,\quad
F=(N+2)^2-4(T-1)^2,
\qquad J_P=\sum_{T,N\in\F_P}\chi(KDF).
\]
\begin{theorem}\label{thm:quadratic-count}
For every odd prime power $P$,
\begin{equation}\label{eq:quadratic-count}
4A_2(P)=4V_P=P^2+1+P\chi(-3)+\chi(-15)-J_P.
\end{equation}
In particular, under the hypotheses defining $C_P$ above,
\[
8b_2=4b_1^2-P^2-1-mP\chi(-3)-m\chi(-15)+mJ_P+4\kappa U_P.
\]
\end{theorem}
\begin{proof}
All sums in this proof are over $(T,N)\in\F_P^2$. We evaluate
\begin{equation}\label{eq:lower-sums}
\begin{gathered}
\sum\chi K=\sum\chi D=\sum\chi F=\sum\chi(KD)=0,\\
\sum\chi(DF)=1,\qquad \sum\chi(KF)=1+P\chi(-3).
\end{gathered}
\end{equation}
The first two sums vanish by summing in $N$. For $F$, put $v=T-1$:
the quadratic sum in $N$ is $-1$ for $v\ne0$ and $P-1$ for $v=0$.
For $KD$, the two linear factors in $N$ coincide exactly at $T=9/2$;
the same quadratic-sum rule proves vanishing, including in
characteristic three.

For the fifth sum, adjoin $z$ with $z^2=D$. Then
$N=(v^2+8-z^2)/4$, and the invertible change
$x=v-z$, $y=v+z$ gives
\[
F=\frac{(x^2-16)(y^2-16)}{16}.
\]
Consequently
\[
\sum(1+\chi D)\chi F
=\left(\sum_x\chi(x^2-16)\right)^2=1.
\]
Since $\sum\chi F=0$, this proves the fifth identity.

For the last sum, instead adjoin $z^2=K$ and put
$x=T-z$, $y=T+z$, so $N=xy/4$. Now
\[
F=\frac{(x-4)(y-4)(xy+4x+4y)}{16}.
\]
Translating $X=x-4$, $Y=y-4$ yields
\[
\sum(1+\chi K)\chi F
=\sum_{X,Y}\chi\bigl(XY((X+8)Y+8X+48)\bigr).
\]
For $X\ne-8,-6$, the sum over $Y$ is $-\chi(X(X+8))$.
At $X=-8$ it is zero; at $X=-6$ its correction to this generic
expression is $P\chi(-12)$. This assertion also holds in characteristic
three, when the correction is zero because $X=-6=0$.
Since $\sum_X\chi(X(X+8))=-1$ and $\chi(-12)=\chi(-3)$,
the last identity in \eqref{eq:lower-sums} follows.

It remains to account for zero discriminants before expanding the
nonsquare indicators. Define
\[
R=\sum(1-\chi K)(1-\chi D)(1-\chi F)-4V_P.
\]
On $D=0$ one has $F=(N-6)^2$. On the two lines $F=0$,
namely $N=2T-4$ and $N=-2T$, one has respectively
$D=(T-5)^2$ and $D=(T+3)^2$.
Thus only the two intersections $(T,N)=(5,6),(-3,6)$ contribute to $R$.
Their contributions are $1-\chi(1)=0$ and $1-\chi(-15)$.
They remain distinct in every odd characteristic. Therefore
$R=1-\chi(-15)$. Expanding the product and using \eqref{eq:lower-sums} gives
\[
4V_P=P^2+(1+P\chi(-3))+1-J_P-(1-\chi(-15)),
\]
which is \eqref{eq:quadratic-count}. The coefficient formula follows from \eqref{eq:second-coefficients} and
$m=4\kappa+1$.
\end{proof}

\subsection{Eliminating the second character sum}
The apparently different sums $J_P$ and $U_P$ are not independent.
The following identity is proved by counting a double cover of an
affine surface in two different sets of coordinates. In particular,
it is an algebraic character-sum identity, not an appeal to a
hypothetical elliptic-curve interpretation of the degree-two count.
\begin{theorem}\label{thm:eliminate-J}
For every finite field of odd order $P$,
\begin{equation}\label{eq:eliminate-j}
J_P=1+P\chi(-3)-U_P.
\end{equation}
Consequently
\begin{equation}\label{eq:a-two}
A_2(P)=\frac{P^2+\chi(-15)+U_P}{4}.
\end{equation}
Under the hypotheses for $S_P,C_P$ above, the quadratic coefficients are
\begin{equation}\label{eq:common-sum}
a_2=\frac{\tau_P^2+U_P}{2},\qquad
b_2=\frac{4b_1^2-P^2-1-U_P+m(1-\chi(-15))}{8}.
\end{equation}
Thus only one character sum, $U_P$, is needed for both quadratic
coefficients. Equivalently, they satisfy the evaluated relation
\begin{equation}\label{eq:linear-relation}
8b_2+2a_2=4b_1^2+\tau_P^2-P^2-1+m(1-\chi(-15)).
\end{equation}
\end{theorem}
\begin{proof}
Put $\epsilon=\chi(-3)$ and
\[
Q(T,N)=12T^2+28(N+1)T+16(N+1)^2+N,
\quad S_0=\sum_{T,N}\chi(NQ(T,N)).
\]
We first evaluate $S_0=P(1+\epsilon)$. If $p\ne3$, the discriminant
of $Q$ as a quadratic in $T$ is $16(N^2-N+1)$. Summing in $T$ and
using $\sum_N\chi N=0$ gives
\[
S_0=P\sum_{N^2-N+1=0}\chi(12N)=P(1+\epsilon).
\]
Indeed $(N+1)^2=3N$ at each such root, so every summand equals one.
If $p=3$, then $Q=(N+1)T+N^2+1$. All inner sums vanish except
the one at $N=-1$, where $NQ=1$ and the contribution is $P$.
This proves the same evaluation in characteristic three.

Define the polynomial
\[
G_0(a,b)=ab(4ab+4a+3b+4)(4ab+3a+4b+4),\qquad
\Sigma_0=\sum_{a,b}\chi(G_0(a,b)).
\]
The change $T=a+b$, $N=ab$ has $1+\chi(T^2-4N)$ ordered preimages.
The polynomial identity
\[
NQ(T,N)\big|_{T=a+b,N=ab}=G_0(a,b)
\]
therefore gives, directly from the definition of $U_P$,
\begin{equation}\label{eq:change-variables}
\Sigma_0=2S_0-U_P=2P(1+\epsilon)-U_P.
\end{equation}

Next count the affine scheme
\[
\mathcal M:\quad x^2=K(T,N),\quad y^2=D(T,N),\quad z^2=F(T,N).
\]
Expansion of its fiber sizes, followed by \eqref{eq:lower-sums}, gives
\[
\#\mathcal M(\F_P)=P^2+2+P\epsilon+J_P.
\]
Conversely, $x,y$ determine $T,N$ uniquely:
\[
T=\frac{x^2-y^2+9}{2},\qquad N=\frac{T^2-x^2}{4}.
\]
Put $u=x-y$, $v=x+y$. Substitution yields the polynomial identity
\begin{align*}
256F(T,N)={}&(u^2-1)(v^2-1)\\
&\cdot\bigl((u-1)(v+1)+16\bigr)
\bigl((u+1)(v-1)+16\bigr).
\end{align*}
Writing the product on the right as $G(u,v)$, we obtain
\[
Z_0:=\sum_{u,v}\chi(G(u,v))=2+P\epsilon+J_P.
\]
On $u,v\ne-1$, make the invertible substitutions
\[
u=\frac{1-a}{1+a},\qquad v=\frac{1-b}{1+b},
\qquad a,b\ne-1.
\]
They satisfy
\[
G(u,v)=\frac{256G_0(a,b)}{(1+a)^4(1+b)^4}.
\]
Every factor removed from the quadratic character is a nonzero square.
The excluded lines $u=-1$ or $v=-1$ contribute zero to $Z_0$.
On the other side, $G_0(-1,b)=b^2$ and $G_0(a,-1)=a^2$; their character
sums are each $P-1$, and their intersection contributes one.
Consequently the exact boundary correction is
\[
Z_0=\Sigma_0-(2P-3).
\]
Together with \eqref{eq:change-variables} this proves \eqref{eq:eliminate-j}, including characteristics three
and five. Substituting \eqref{eq:eliminate-j} into \eqref{eq:quadratic-count} proves \eqref{eq:a-two}, and substitution into
\eqref{eq:second-coefficients} gives \eqref{eq:common-sum} and \eqref{eq:linear-relation}.
\end{proof}

\begin{remark}[Elliptic factors and point counts]
At $P=3$, the factor $1+u+3u^2$ in $S_3=N$ is the Frobenius
polynomial of the elliptic curve $y^2=x^3+2x^2+1$ over $\F_3$,
which has five rational points. By semisimplicity and Tate's
isogeny theorem for abelian varieties over finite fields
\cite{Tate}, this elliptic curve is an isogeny factor of
$\operatorname{Jac}(\mathcal B_3)$ over $\F_3$.

For $p\ne3,5$, the equality $\#\mathcal B_P(\F_P)=\#E(\F_P)$
compares the first Frobenius traces over $\F_P$. A comparison over
extension fields would require higher Frobenius power traces.
The quadratic-coefficient identities above express the next
coefficients of $S_P$ and $C_P$ in terms of the finite-field sum $U_P$.
\end{remark}

\section{A quotient of covering-curve zeta functions}\label{sec:quotient}
\begin{remark}
Let $X$ be the symplectic frame torsor of $\mathcal V$, let
$Z=\{\pm I\}$, and let $X_Z$ and $X_T$ be the smooth projective
compactifications of $X/Z$ and $X/T$, respectively. The latter
quotient uses the Coxeter torus $T$ lifting $T_0$.
For each boundary point of $X_T$, the decomposition group of
$X_Z\to X_T$ is the image in $T/Z$ of the intersection of $T$
with a conjugate of a decomposition group of $X$ over infinity.
Each such conjugate preserves a Lagrangian plane: the original
group preserves the positive-valuation plane. Every element of
$T\setminus Z$ acts irreducibly, so this intersection is contained
in $Z$. Thus every boundary point splits completely in the
degree-$m$ map $X_Z\to X_T$.

On the affine line, Frobenius acts on the fibres of $X/T$ and $X/Z$
with permutation characters $\operatorname{Ind}_{T_0}^S1$ and
$\operatorname{Reg}_S$, respectively. Summing their difference
$m\operatorname{Ind}_{T_0}^S1-\operatorname{Reg}_S=4m1_{\rm irr}$
over $\F_{P^H}$ gives the affine point-count difference. Complete
splitting makes the corresponding boundary contributions cancel,
so
\[
4m A_H(P)=m\#X_T(\F_{P^H})-\#X_Z(\F_{P^H}),\qquad
\mathcal R_P(u)=\frac{Z(X_T,u)^m}{Z(X_Z,u)}.
\]
The second identity follows by exponentiating the first with
weights $u^H/H$.
\end{remark}

\appendix
\section{Higher-rank symplectic flag covers and Steinberg factors}
\label{app:higher-flags}

The flag curve of Section~\ref{sec:flag} is an instance of a construction
which applies to symplectic local systems of arbitrary even rank.  The
underlying reason for considering full flags is the Solomon--Tits
realization of the Steinberg representation in the top reduced homology
of the spherical building.  We record the general construction and its
cohomological consequence, without assuming that the required
higher-rank local systems arise from the skew polynomials considered
in this paper.

Let $k=\F_P$, let $X$ be a smooth geometrically connected affine curve
over $k$, and let $\mathcal V$ be a lisse sheaf of $2n$-dimensional
$\F_P$-vector spaces on $X$, equipped with a nondegenerate alternating
pairing.  Assume that its geometric monodromy is
$G=\Sp_{2n}(\F_P)$; in particular the arithmetic monodromy is also $G$.
Fix a symplectic fibre $V$.  Write $\Delta$ for the set of $n$ simple
roots of $G$, and for $J\subseteq\Delta$ let $P_J$ be a standard
parabolic subgroup, with $P_\varnothing=B$ a Borel subgroup and
$P_\Delta=G$.  The character $\operatorname{Ind}_{P_J}^G\mathbf1$ is the permutation
character of $G$ on the coset set $G/P_J$. The finite
transitive $G$-set $G/P_J$ consists of
isotropic partial flags of the corresponding type.  By the monodromy
action it defines a finite \'etale cover
\[
 \pi_J:Y_J\longrightarrow X.
\]
Every $Y_J$ is geometrically connected.  The cover $Y_\varnothing$
parametrizes complete isotropic flags
\[
 0=V_0\subset V_1\subset\cdots\subset V_n\subset V,
 \qquad \dim_{\F_P}V_i=i,\quad V_n=V_n^\perp,
\]
and $Y_\Delta=X$.  Although the flags have length $n$, each $Y_J$
remains a curve, since it is finite over $X$.

\begin{proposition}[The building identity for flag covers]
\label{prop:higher-building}
Let $\mathrm{St}_G$ be the Steinberg character of $G$, and let
$F_{x,h}\in G$ denote the Frobenius conjugacy class of the fibre at
$x\in X(\F_{P^h})$.  Then
\begin{equation}\label{eq:higher-building-character}
 \mathrm{St}_G=\sum_{J\subseteq\Delta}(-1)^{|J|}
       \operatorname{Ind}_{P_J}^{G}\mathbf 1,
\end{equation}
and, for every $h\ge1$,
\begin{equation}\label{eq:higher-building-point-count}
 \sum_{x\in X(\F_{P^h})}\mathrm{St}_G(F_{x,h})
 =\sum_{J\subseteq\Delta}(-1)^{|J|}
       \#Y_J(\F_{P^h}).
\end{equation}
If $\mathcal E_{\mathrm{St}}$ is the lisse
$\overline{\mathbb Q}_\ell$-sheaf associated with $\mathrm{St}_G$
($\ell\ne p$), then its $L$-function satisfies
\begin{equation}\label{eq:higher-building-L}
 L(X,\mathcal E_{\mathrm{St}},u)
 =\prod_{J\subseteq\Delta}Z(Y_J,u)^{(-1)^{|J|}}.
\end{equation}
\end{proposition}
\begin{proof}
The simplices of the spherical building of $G$ are the proper
isotropic flags; its augmented simplicial chain complex has, in each
degree, the permutation modules on flags of the relevant types.
The Solomon--Tits theorem (see, e.g., \cite{BrownBuildings}) identifies its sole nonzero reduced
homology with $\mathrm{St}_G$, giving
\eqref{eq:higher-building-character} by the Euler characteristic.
The value of $\operatorname{Ind}_{P_J}^{G}\mathbf1$ at $F_{x,h}$ counts
$F_{x,h}$-fixed flags in $G/P_J$, equivalently the
$\F_{P^h}$-points of the fibre of $Y_J$ above $x$.  Summing over
$x$ gives \eqref{eq:higher-building-point-count}.  Finally,
apply the exponential definition of the $L$-function to the
trace identity for all $h$; the zeta function of $Y_J$ is the
$L$-function of the permutation sheaf $(\pi_J)_*\overline{\mathbb Q}_\ell$.
\end{proof}

For comparison with projective curves, write $\overline Y_J$ for the
smooth projective compactification of $Y_J$ and
$D_J=\overline Y_J\setminus Y_J$ for its reduced geometric boundary.
Let
\[
 P_J^{\mathrm{geom}}(u)
 =\det(1-uF\mid H^1(\overline Y_{J,\bar k},\overline{\mathbb Q}_\ell)),
 \qquad
 B_J(u)=\det(1-uF\mid H^0(D_{J,\bar k},\overline{\mathbb Q}_\ell)).
\]
The finite set $D_J$ carries its natural Frobenius permutation action.
The localization sequence gives
\[
 \det(1-uF\mid H_c^1(Y_{J,\bar k}))
 =P_J^{\mathrm{geom}}(u)\,\frac{B_J(u)}{1-u}.
\]
Since $Z(Y_J,u)=Z(\overline Y_J,u)/Z(D_J,u)$ and
$\sum_J(-1)^{|J|}=0$, equation~\eqref{eq:higher-building-L}
becomes
\begin{equation}\label{eq:higher-projective}
 L(X,\mathcal E_{\mathrm{St}},u)
 =\prod_{J\subseteq\Delta}
 \left(P_J^{\mathrm{geom}}(u)\,B_J(u)\right)^{(-1)^{|J|}}.
\end{equation}
In particular, the boundary correction is indispensable in general.
If $X=\mathbb A^1$, the absence of compactly supported cohomology
in degrees $0$ and $2$ for the nontrivial Steinberg local system gives
\[
 L(X,\mathcal E_{\mathrm{St}},u)
 =\det(1-uF\mid H_c^1(X_{\bar k},\mathcal E_{\mathrm{St}})).
\]

\begin{remark}[The rank-two specialization and a higher-rank question]
\label{rem:higher-rank-question}
When $n=2$, the four terms of~\eqref{eq:higher-building-character}
are, respectively, full flags, isotropic lines, Lagrangian planes,
and the trivial representation:
\[
 \mathrm{St}_G
 =\mathbf1_{\mathrm{Flags}}-\mathbf1_{\mathrm{Lines}}
  -\mathbf1_{\mathrm{Lag}}+\mathbf1.
\]
In the specific family studied in Section~\ref{sec:flag}, the
lower-dimensional flag covers have elementary cohomological factors,
and their boundary contributions cancel the extraneous factors of the
full-flag cover.  This is why the projective full-flag curve
$\mathcal B_P$ alone realizes $L_{\mathrm{St}}$.
For $n\ge3$, the building identity still produces a canonical
\emph{alternating product} of zeta functions, but does not imply
that a single smooth projective flag cover has zeta numerator
$L(X,\mathcal E_{\mathrm{St}},u)$.  Obtaining explicit plane models,
computing inertia orbits and genera, and determining when such a
single-curve realization occurs are separate problems.  No assertion
of higher-rank symplectic monodromy for the skew-polynomial family is
made here.
\end{remark}

\section*{Declaration of generative AI assistance}
OpenAI ChatGPT (GPT-6) and Codex were used to assist with mathematical
exploration, reviewing arguments, checking calculations and references,
preparing and running finite-field verification programs, and improving
the exposition and language of this manuscript. The authors independently
reviewed the AI-assisted material and take full responsibility for the
correctness of the mathematical statements, proofs, computational results,
references, and final manuscript.



\begin{thebibliography}{20}
\bibitem{AS} A. Abbes and T. Saito,
\emph{Local Fourier transform and epsilon factors},
Compos. Math. \textbf{146} (2010), no.~6, 1507--1551.
\url{https://doi.org/10.1112/S0010437X09004631}.
\bibitem{Breeding} J. Breeding II,
\emph{Irreducible characters of $\mathrm{GSp}(4,q)$ and dimensions of spaces of fixed vectors},
Ramanujan J. \textbf{36} (2015), no.~3, 305--354.
\url{https://doi.org/10.1007/s11139-014-9622-3}.
\bibitem{BrownBuildings} K. S. Brown,
\emph{Buildings}, Springer, New York, 1989, Chapter IV.
\bibitem{CL} X. Caruso and J. Le Borgne,
\emph{A new faster algorithm for factoring skew polynomials over finite fields},
J. Symbolic Comput. \textbf{79} (2017), part~2, 411--443.
\url{https://doi.org/10.1016/j.jsc.2016.02.016}.
\bibitem{DL} P. Deligne and G. Lusztig,
\emph{Representations of reductive groups over finite fields},
Ann. of Math. (2) \textbf{103} (1976), no.~1, 103--161.
\url{https://doi.org/10.2307/1971021}.
\bibitem{Geck} M. Geck,
\emph{Some applications of CHEVIE to the theory of algebraic groups},
Carpathian J. Math. \textbf{27} (2011), no.~1, 64--94.
\url{https://arxiv.org/abs/1008.1464}.
\bibitem{Generic} GenericCharacterTables.jl, generic unipotent table
\texttt{uniCSp4.1}, based on Shinoda's character table \cite{Shinoda},
\url{https://github.com/oscar-system/GenericCharacterTables.jl/blob/master/data/Tables/C2/uniCSp4.1.jl}.
\bibitem{HTtrinomials} T. Hayata and T. Tsushima,
\emph{Counting irreducible skew polynomials $x^{2n}+ax^n+1$},
arXiv:2611.04767 (2026).
\url{https://arxiv.org/abs/2611.04767}.
\bibitem{GeckMalle} M. Geck and G. Malle,
\emph{The Character Theory of Finite Groups of Lie Type: A Guided Tour},
Cambridge Studies in Advanced Mathematics 187, Cambridge University Press, 2020.
\bibitem{KT} N. M. Katz and P. H. Tiep,
\emph{Rigid local systems and finite symplectic groups},
Finite Fields Appl. \textbf{59} (2019), 134--174.
\url{https://doi.org/10.1016/j.ffa.2019.05.008}.
\bibitem{Laumon} G. Laumon,
\emph{Transformation de Fourier, constantes d'\'equations fonctionnelles
et conjecture de Weil},
Publ. Math. Inst. Hautes \`Etudes Sci. \textbf{65} (1987), 131--210.
\url{https://doi.org/10.1007/BF02698937}.
\bibitem{LMFDB15} The LMFDB collaboration, elliptic-curve isogeny class
$15.a$ and newform $15.2.a.a$,
\url{https://www.lmfdb.org/EllipticCurve/Q/15/a/},
\url{https://www.lmfdb.org/ModularForm/GL2/Q/holomorphic/15/2/a/a/}.
\bibitem{Ore} O. Ore,
\emph{Theory of non-commutative polynomials},
Ann. of Math. (2) \textbf{34} (1933), no.~3, 480--508.
\url{https://doi.org/10.2307/1968173}.
\bibitem{Shinoda} K.-i. Shinoda,
\emph{The characters of the finite conformal symplectic group,
$\operatorname{CSp}(4,q)$},
Comm. Algebra \textbf{10} (1982), 1369--1419.
\url{https://doi.org/10.1080/00927878208822782}.
\bibitem{Tate} J. Tate,
\emph{Endomorphisms of abelian varieties over finite fields},
Invent. Math. \textbf{2} (1966), 134--144.
\bibitem{Srinivasan} B. Srinivasan,
\emph{The characters of the finite symplectic group $\operatorname{Sp}(4,q)$},
Trans. Amer. Math. Soc. \textbf{131} (1968), 488--525.
\url{https://doi.org/10.2307/1994960}.
\end{thebibliography}
\end{document}